\documentclass{ajour}
\usepackage{amsfonts}
\usepackage{amsmath}
\usepackage{amssymb}
\input{epsf}

\newcommand{\la}{\langle}
\newcommand{\ra}{\rangle}

\newcommand{\Be}{{\beta}}
\newcommand{\al}{{\alpha}}
\newcommand{\im}{\mathop{\rm{ im\, }}}
\newcommand{\lowim}{{\mathop{\rm{Lowim}}}}
\usepackage{verbatim}
\newcommand*{\clo}[1]{\overline{#1}}

\begin{document}

\title{On the topology of components \\ of some
Springer fibers and their relation to Kazhdan-Lusztig theory}
\thanks{This paper represents part of a dissertation written under the 
direction of Prof. Bob MacPherson at the Institute for Advanced Study. 
The author would like to take this opportunity to express gratitude for his 
unwavering support, guidance and enthusiasm. The author would also like to
acknowledge Tom Braden, Misha Grinberg, Mark Goresky, Arun Ram, and the 
anonymous referee for many helpful discussions and comments. The author was 
partially supported by an NDSEG fellowship during his graduate studies at
Princeton University.}
\date{April 17, 2002}

\author{Francis Y. C. Fung}
\authorrunninghead{Francis Y. C. Fung}
\titlerunninghead{Topology of Some Springer Fibers}

\email{fycfung@alumni.princeton.edu}


\abstract{We describe the irreducible components of 
Springer fibers for hook and two-row nilpotent elements of 
$
\mathfrak{gl}_n(\mathbb{C})
$
as iterated bundles of flag manifolds and Grassmannians. We then
relate the topology (in particular, the intersection homology
Poincar\'{e} polynomials) of the intersections of these components with
the inner products of the Kazhdan-Lusztig basis elements of irreducible
representations of the rational Iwahori-Hecke algebra of type $A$
corresponding to the hook and two-row Young shapes. This work has been submitted to Advances in Mathematics (Academic Press) for possible publication. Copyright may be transferred without notice, after which this version may no longer be accessible.
}



\section{Introduction}
Let $V$ be a finite-dimensional complex vector space. 
A nilpotent linear map $N:V \to V$ is said to fix a flag $F =\{ F_0
\subset  F_1 \subset \dots \subset F_{n-1} \subset V\}$ if $NF_i \subseteq F_{i-1}$
for each $i$.
 The variety $\mathcal{B}_N$ of all flags in 
the flag manifold $Fl(V)$ fixed by a nilpotent map $N$ is a
 \emph{Springer fiber}. Such varieties arise as fibers of Springer's resolution of 
singularities of the nilpotent cone of a reductive algebraic group $G$.

Springer\cite{Spr2} discovered a method of constructing irreducible
representations of Weyl groups on the top homology of
$\mathcal{B}_N$. For each irreducible representation, his
construction yields a distinguished basis given by homology classes of
the components of the Springer fiber $\mathcal{B}_N$. 
However, since their mere existence yields the distinguished
basis, it seems that efforts to understand them have not focused on
computations of their internal topological structure or that of their
intersections.  Only a few papers
(such as Spaltenstein\cite{Spa}, Vargas\cite{V}, Wolper\cite{Wo},
Lorist\cite{Lo}, G\"uemes\cite{Gu}) have studied the topology of the
components of the Springer fibers $\mathcal{B}_N$ and their
intersections. We extend some of these results to
describe the homological structure of components of $\mathcal{B}_N$
and their pairwise intersections for certain types of nilpotent maps
$N$ (those corresponding to hook and two-row shape partitions). In these cases, the components are
nonsingular, and in fact are iterated bundles of flag manifolds and
Grassmannians. For more general nilpotent maps $N$, the components can
be singular (see Vargas\cite{V} and Spaltenstein\cite{Spa}) and much
more complicated.

We also relate our computations to the structure of the 
Kazhdan-Lusztig bases of certain representations
of Iwahori-Hecke algebras of type $A$. The inner
products of these basis vectors, suitably normalized, are
polynomials in $t$ and $t^{-1}$ that are invariant under the map
$t \to t^{-1}$. We show that for irreducible representations labeled by a
hook or two-row shape, the (suitably normalized) inner products equal
the intersection homology Poincar\'e polynomials of pairwise
intersections of irreducible components of Springer fibers of the
general linear group. We believe it would be very
interesting to understand how our results might generalize to bases of
other Kazhdan-Lusztig representations of type $A$.

\section{Some properties of nilpotent maps and Springer fibers} 
We record some properties of nilpotent maps and of the space of all flags
$\mathcal{B}_N$ fixed by a nilpotent map $N$, which is called the
\emph{Springer fiber} of $N$.  A theorem of Vargas and Spaltenstein
\cite{V} \cite{Spa}
decomposes the space
$\mathcal{B}_N$ into a disjoint union of locally closed subspaces,
whose closures are the irreducible components of the space
$\mathcal{B}_N$.

Let  $N:V \to V$ be a nilpotent map of a vector space $V$ over
$\mathbb{C}$. Let $b$ be the least positive integer for which $N^b =
0$. Then we have two filtrations of subspaces on $V$:
the image filtration $ \im N^b =0 \subset \im N^{b-1} \subset \im
N^{b-2}\subset \dots \subset
\im N^{1} \subset V=\im N^0$ and the kernel filtration $\ker N^0 = 0 \subset \ker N
\subset\ker N^2 \subset \dots \subset \ker N^{b-1} \subset V= \ker
N^b$  (with proper inclusions).

\begin{lemma}
\label{kerlem} For a nilpotent map $N$, we have $N^{-1}(\im N^k) =\ker
N + \im N^{k-1}.$
\end{lemma}

\begin{proof} If $N(v) \in \im N^k$ then $N(v) = N^k(w)$ so $N(v-
N^{k-1}w) =0$.
\end{proof}

Note that if $\im N^{k-1}$ contains $\ker N$ then $N^{-1}(\im
N^{k}) = \im N^{k-1}$; otherwise $N^{-1}(\im N^k)$ is strictly larger.
Also note that $N(\ker N^{i+1}) \subseteq \im N \cap \ker N^{i}$.

\begin{definition} Let $F_i$ be a subspace of $V$ that is taken into
itself by the map $N$. Then there is a map $N_i:V/F_i \to V/F_i$
induced by $N$. We call the map $N_i$ a \emph{quotient map} of $N$.
\end{definition}

\begin{lemma}
\label{quotlem} The image $\im N_i$ of the quotient map $N_i:V/F_i \to
V/F_i$ is equal to $(\im N + F_i) / F_i$.  Similarly, $(\im N^k
+F_i)/F_i= \im N_i^k$.
\end{lemma}

\begin{proof}
If $N(v) \in \im N$ then $N(v) +F_i \in \im N_i$. On the other hand,
if $w + F_i \in \im N_i$ then the coset $w +F_i$ equals the coset
$N(v) +F_i$ for some $v \in V$, so $w + F_i$ is clearly in $\im N+
F_i$. Then a subspace of $V$ that contains $F_i$ and whose projection
to $V/F_i$ is $\im N_i$ must be $\im N + F_i$. The same holds true for
the nilpotent map $N^k$.
\end{proof}

\begin{lemma}
\label{invkerlem} The kernel of the quotient map $\ker N_i$ equals
$N^{-1}(F_i)$.
\end{lemma}

\begin{proof} The kernel $\ker N_i$ is given by those vectors whose
image under $N_i$ is $0+F_i$, which is exactly $N^{-1}(F_i)$.
\end{proof}

\begin{lemma} We have the containment $\ker N \supseteq \im N^{b-1}$,
but $\ker N \not\supseteq
\im N^{b-j}$ for $j > 1$.
\end{lemma}
\begin{proof} This is obvious from definition of $b$.
\end{proof}

\begin{lemma} If $j$ is the largest integer for which $\im N^{j}$ is
not contained in $F_i$, then $j$ is the largest integer for which
$(\im N^{j}+ F_i)/ F_i$ is a nonzero image of the quotient map $N_i$.
\end{lemma}
\begin{proof} Suppose $v \in \im N^j$ but is not in $F_i$. Then $v +
F_i$ is not zero in $V/F_i$, so $(\im N^{j}+ F_i)/ F_i$ is a nonzero
image of $N_i$. Similarly, if $v + F_i$ is a nonzero element of $(\im
N^{j}+ F_i)/ F_i$, then there exists an element $0 \neq w \in V$ with
$w \in \im N^j$ and $w \in v + F_i$. \end{proof}

Now we discuss the Springer fiber $\mathcal{B}_N$ of $N$, which is the
variety of flags fixed by the nilpotent element $N$.  The ranks of the
Jordan blocks of the nilpotent map $N$ determine a partition of
$n$. We form a Young shape from this partition by using this partition
as the lengths of the rows (opposite to the
conventions of Vargas\cite{V}).  Let the number of columns of the Young
shape be $b$; then $N^b = 0$ and $N^{b-1} \neq 0$.

\begin{definition} Given a flag $F$ with subspaces $\{0\}=F_0
\subseteq F_1 \subseteq \dots
\subseteq F_n =V$, we say that
\emph{$N$ fixes $F$} if $NF_i$ is contained in $F_{i-1}$ for all $i$.
\end{definition}

\begin{definition} We denote by $\mathcal{B}_N$ the set of all flags
fixed by $N$, and call it the Springer fiber of $N$. It is an
algebraic subvariety of the flag manifold $\mathcal{B}$.
\end{definition}

Recall that \emph{Young shape} on $n$ boxes is a collection
of $n$ boxes arranged in left justified rows of lengths $t_1 \geq
\dots \geq t_k$. A \emph{standard tableau} (or \emph{Young tableau})
 on a Young shape $\tau$ is constructed by filling in the $n$ boxes
 with the numbers $1, \dots, n$ such that the numbers are decreasing
 from left to right in each row, and decreasing from top to bottom in
 each column. Note that many authors use increasing rows and
 columns. We typically use $A$ and $B$ to denote standard Young
 tableaux.  Denote by $A_i$ the tableau obtained by deleting
the numbers $1, \dots, i$ in the tableau $A$ and subtracting $i$ from
the remaining numbers.

The following theorem of Vargas and Spaltenstein
gives a decomposition of the Springer fiber $\mathcal{B}_N$ into a
disjoint union of locally closed subsets, whose closures comprise the
irreducible components of $\mathcal{B}_N$.

\begin{theorem}\label{SVthm} (Vargas\cite{V}, Spaltenstein\cite{Spa})
Let $N$ be a nilpotent map. Then given a standard tableau $A$ on the
Young shape of $N$, we construct a locally closed subset $SV(A)$ of
the Springer fiber $\mathcal{B}_N$, whose closure $\clo{SV(A)}$ is an
irreducible component of $\mathcal{B}_N$.  We have a partition
$\mathcal{B}_N = \cup_A SV(A)$ of the Springer fiber into disjoint
locally closed subsets.  Thus the number of irreducible components of
$\mathcal{B}_N$ is equal to the number of standard tableaux on the
Young shape of $N$. In addition, the components are all of the same
dimension. In fact, if the lengths of the columns of the Young shape
of $N$ are $n_1, n_2, \dots, n_b$, then the dimension of each
component is
\[
\sum_i \frac{n_i(n_i-1)}{2}.
\]
\end{theorem}

\begin{proof} Suppose $A$ is a Young tableau on the Young shape of
$N$. Then we inductively 
specify a subset of $\mathcal{B}_N$ corresponding to $A$,
which we denote $SV(A)$ (for Spaltenstein-Vargas),
\label{SVdef} by describing how to choose $F_1$, then
$F_2/F_1$, and so forth. A flag $F$ is in the the subset $SV(A)$ if
each subspace of $F$ satisfies the following conditions.

Suppose the number $i$ appears in the $c(i)$-th column in $A$. Then
(recalling that $F_0 = 0$) the first subspace $F_1$ must satisfy $F_1
\subset (N^{-1}(F_0) \cap (\im N^{c(1) -1} - \im N^{c(1)})$. In other
words, $F_1$ must be in the kernel of $N$, and it must be in the
($c(1)- 1$)st image of $N$ but not in any higher image.

Then, for any $F_1$ satisfying the above condition, the induced map
$N_1:V/F_1 \to V/F_1$ will have the same Young shape; this shape is
the shape obtained by deleting the number $1$ in the tableau $A$.

We now choose $F_2/F_1 \subset V/F_1$, using the above procedure with
$N_1$ in place of $N$ and $A_1$ in place of $A$. Note that
$N_1^{-1}(0+F_1) = N^{-1}(F_1)/F_1$. We continue inductively, choosing
$F_{i+1}/F_{i}$ so that $F_{i+1}/F_i \subset (N^{-1}(F_i)/F_i) \cap
(\im N_{i}^{c(i) - 1} - \im N_{i}^{c(i)})$. Note that any such choice
of $F_{i+1}$ yields a quotient map $N_{i+1}$ with Young shape
$A_{i+1}$.

See Vargas\cite{V} and Spaltenstein\cite{Spa}
for the proof that this constructs a locally closed subset
of $\mathcal{B}_N$ with the properties claimed in the theorem. The
proof is essentially an explicit calculation with the Jordan form of
$N$.

Vargas(\cite{V}, Proposition 2.2) shows that this set $SV(A)$ is
exactly the set of flags such that
\begin{align*} F_i \subset N^{-1}(F_{i-1})\\ F_{i} \subset F_{i-1} +
\im N^{c(i) -1}.
\end{align*}

\end{proof}

\section{Determination of the topology of the irreducible components
of Springer fibers for nilpotent maps of hook type for
$GL_n(\mathbb{C})$.}

Suppose $N$ is a nilpotent map $V \to V$ whose Jordan form has at most
one Jordan block of rank $> 1$. Then $N$ is said to be of hook type.
For hook type nilpotent maps $N$, we can characterize the components
of the Springer fiber $\mathcal{B}_N$ entirely in terms of the image
and kernel filtrations of $V$. The components and their intersections
will turn out to be nonsingular. In fact their homology Poincar\'e
polynomials factor as products of Poincar\'e polynomials of
Grassmannians and flag manifolds.

We will describe each component of the Springer fiber of a nilpotent
map $N$ of hook type by expressing the component as a sequence of
fiber bundles with progressively simpler bases and fibers.

\begin{definition} A space $X_1$ is an \emph{iterated fiber bundle} of
base type $(B_1, \dots, B_n)$ if there exist spaces $X_1, B_1, X_2,
B_2, \dots, X_n, B_n, X_{n+1}=pt$ and maps $p_1, p_2,
\dots, p_{n}$ such that $p_j:X_j \to B_j$ is a fiber bundle with
typical fiber $X_{j+1}$.
\end{definition}

The following two lemmas are straightforward.

\begin{lemma}: The flag manifold $Fl(V)$ admits a map to the
Grassmannian $G_i(V)$ via $F \mapsto F_i$. The fiber of this map is a
product $Fl(F_i) \times Fl(V/F_i)$.
\end{lemma}

\begin{lemma}\label{IKlem} Consider the variety of flags $X(I, K)$ in
an $n$-dimensional vector space $V$ such that $F_i$ contains an
$a$-dimensional subspace $I$ and is contained in a $b$-dimensional
subspace $K$. This variety $X(I,K)$ admits a map $p: X(I, K) \to
G_{i-a}(K/I)$ via $F \mapsto F_i/I \subset K/I$ that makes $X(I,K)$
the total space of a fiber bundle.  The typical fiber of the map $p$
is a product $Fl(F_i) \times Fl(V/F_i)$. In particular, the variety
$X(I,K)$ is nonsingular.
\end{lemma}

We define some notation to  simplify our intersection homology
computations. 
Define $[n]$ by \[[n] :=
t^{-(n-1)}(1+ t^2+ t^4 + \dots +t^{2(n-1)}).\] Then we define \[[n]!
:= [1][2]\dots[n]
\hbox{ and } \binom{[n]}{[k]} := \frac{[n]!}{ [k]![n-k]!}.\]

 The polynomial $[n]$ is essentially the $t$-analogue of
the number $n$, but shifted to be symmetric around the degree $0$
term.  Thus the flag manifold $Fl(V)$ has intersection homology
Poincar\'e polynomial $[n]!$ and $G_k(V)$ has intersection homology
Poincar\'e polynomial $\binom{[n]}{[k]}$. 

\begin{corollary} Let $X(I)$ be the variety of flags in $Fl(V)$ such
that the the subspace $F_i$ contains an $a$-dimensional subspace
$I$. Then $X(I)$ has intersection homology Poincar\'e polynomial
$\binom{[n-a]}{[i-a]} [i]! [n-i]!$.
\end{corollary}
\begin{proof} Clear from Lemma \ref{IKlem}, since complex flag
manifolds and complex Grassmannians have only even-dimensional
homology so the Leray-Serre spectral sequence for $p: X(I) \to G_{i-a}(V/I)$
with fiber $Fl(F_i) \times Fl(V/F_i)$ collapses.
\end{proof}

\begin{lemma} \label{hookimlem} Let $N$ be a nilpotent map of hook
type with Jordan blocks of size $(b, 1, \dots, 1)$. Then for all $0 <
i < b$, we have $\im N^{b-i} = \ker N^{i}\, \cap \, \im N$, which
implies $N(\ker N^{i+1}) \subset \im N^{b-i}$.
\end{lemma}
\begin{proof} This follows from inspection of the Jordan form of $N$.
\end{proof}

We now decompose  each component of the
Springer fiber $\mathcal{B}_N$ for a hook type nilpotent map $N$ as an
iterated bundle with nonsingular bases and fibers.  Recall that the
number $b$ is the least positive integer with $N^b = 0$. 
Let $b(i)$ be the least positive integer with $N_i^{b(i)} = 0$.

\begin{theorem} \label{hookchar} Suppose we are given a nilpotent map
$N$ of hook type and a Young tableau $A$ on the Young shape of $N$
with $n, i_{b-1}, \dots, i_1$ on the top row (where by convention $i_b
= n$ and $i_0 = 0$). Then the component $K_A$ of the Springer fiber
$\mathcal{B}_N$ is an iterated bundle with $B_{2j-1} = G_{i_{j}-
i_{j-1}-1}(\ker N_{i_j}/\im N^{b(i_j)-1}_{i_j})$ and $B_{2j} =
Fl(F_{i_j}/F_{i_{j-1}})$, where $j = 1, 2, \dots, b-1$, and $B_{2b-1}$
is a full flag manifold $Fl(V/F_{i_{b-1}})$.
\end{theorem}

The proof will be broken up into a series of lemmas and propositions.

\begin{proposition} \label{vargasprop} (Vargas\cite{V}) Suppose we are
given a nilpotent map $N$ of hook type and a Young tableau $A$ on the
Young shape of $N$ with $n, i_{b-1}, \dots, i_1$ on the top row.  Then
the component $K_A = \clo{SV(A)}$ of the Springer fiber
$\mathcal{B}_N$ consists of all flags $F$ in $Fl(V)$ such that
\begin{gather*}
\im N^{b-1} \subseteq F_{i_1} \subseteq \ker N \\
\im N^{b-2} \subseteq F_{i_2} \subseteq \ker N^2\\
\im N^{b-3} \subseteq F_{i_3} \subseteq \ker N^3\\
\dots\\
\im N^{1} \subseteq F_{i_{b-1}} \subseteq \ker N^{b-1}.
\end{gather*}
\end{proposition}
\begin{proof} This is Vargas\cite{V}, Theorem 4.1. The proof consists
of an explicit limiting argument using the structure of the
Spaltenstein-Vargas subset $SV(A)$.
\end{proof}

\begin{lemma} \label{quotimkerlem} Let $N$ be a nilpotent map of hook
type, and let $F_{i_1}$ be a subspace of $V$ with $NF_{i_1} \subset
F_{i_1}$ and $\im N^{b-1} \subset F_{i_1} \subset \ker N$. Then we
have
\begin{gather*}
\ker N_{i_1}^d = \ker N^{d+1}/F_{i_1} \\
\im N_{i_1}^d = (\im N^d + F_{i_1})/ F_{i_1}.
\end{gather*}
\end{lemma}

\begin{proof} We describe $\ker N_{i_1}^d \subseteq V/F_i$ as
follows. If $v \in \ker N^{d+1}$ then $N^{d}(v)
\in \im N \cap \ker N = \im N^{b-1}$. Since $\im N^{b-1} \subseteq
F_{i_1}$, we have $N^{d}(v)
\in F_{i_1}$, so $v +F_{i_1} \in \ker N_{i_1}^d$. On the other hand,
if $N_{i_1}^d(v + F_{i_1}) = 0 +F_{i_1}$ then $N^d(v) \in
F_{i_1}$. Now by Lemma \ref{hookimlem}, the subspace $F_{i_1}$
contains $\im N^{b-1}$ but no other element of $\im N$. Thus
$N^{d+1}(v) =0$, and so $v \in
\ker N^{d+1}$. Thus if a subspace $W \subset V$ contains $F_{i_1}$
then $W/F_{i_1} \subseteq \ker N_{i_1}^{d}$ iff $W \subseteq \ker
N^{d+1}$.

By Lemma \ref{quotlem} $\im (N^d +F_{i_1}) / F_{i_1} = \im
N_{i_1}^d$. Thus if $W \subset V$ contains $F_{i_1}$, then $W/F_{i_1}$
contains $\im N_{i_1}^d$ iff $W$ contains $\im N^d$.
\end{proof}

\begin{lemma} \label{smallerspringerlem} Suppose we are given a
nilpotent $N$ of hook type and a Young tableau $A$ on the Young shape
of $N$ with $n, i_{b-1}, \dots, i_1$ on the top row (where by
convention $i_b = n$ and $i_0 = 0$).  Then the component $K_A$ of the
Springer fiber $\mathcal{B}_N$ admits a map $p_1$ to the Grassmannian
$G_{i_1-1}(\ker N / \im N^{b-1} )$.  The fiber $X_2$ of the map
$p_1:K_A \to G_{i_1-1}(\ker N /
\im N^{b-1} )$ admits a map $p_2: X_2 \to Fl(F_{i_1})$. The fiber of
$p_2$ can be identified with a component of a Springer fiber of the
quotient map $N_{i_1}:V/F_{i_1} \to V/F_{i_1}$, where the component is
associated to the standard tableau $A_{i_1}$.
\end{lemma}

\begin{proof} The existence of the map $p_1$ follows from Lemma
\ref{IKlem} with $I = \im N^{b-1}$ (which is a $1$-dimensional space)
and $K = \ker N$ (which is an $n-b+1$-dimensional space containing
$\im N^{b-1}$). Let $B_1$ be the Grassmannian $G_{i_1-1}(\ker N/\im
N^{b-1})$. The fiber $X_2$ of the map $p_1: K_A \to B_1$ consists of
all flags in the component with a fixed subspace $F_{i_1}$. We define
the map $p_2:X_2 \to Fl(F_{i_1})$ by taking $F \in X_2$ and forgetting
all subspaces of $F$ larger than $F_{i_1}$. By inspecting Proposition
\ref{vargasprop} we see that $p_2$ is surjective and indeed 
a fiber bundle projection.

The fiber $X_3$ of this map $p_2$ is the set of all flags in the
component $K_A$ with fixed subspaces $F_1, \dots, F_{i_1}$. Then $X_3$
maps bijectively to a subset of $Fl(V/F_{i_1})$ via the map $F \to
F'$, where $F_j' = F_{j +i_1}/F_{i_1}$. We now show that $X_3$ is the
component $K_{A_{i_1}}$ of the Springer fiber of the quotient map
$N_{i_1}$ on $V/F_{i_1}$, by showing that $X_3$ satisfies the
characterization of Proposition \ref{vargasprop}.

By Lemma \ref{quotimkerlem}, the fiber $X_3$ of $p_2$ is in bijection
with the set of flags $F'
\in Fl(V/F_{i_1})$ such that \begin{gather*}
\im N_{i_1}^{b-2} \subseteq F_{i_2-i_1}' \subseteq \ker N_{i_1} \\
\im N_{i_1}^{b-3} \subseteq F_{i_3-i_1}' \subseteq \ker N_{i_1}^2\\
\im N_{i_1}^{b-4} \subseteq F_{i_3-i_1}' \subseteq \ker N_{i_1}^3\\
\dots\\
\im N_{i_1}^{1} \subseteq F_{i_{b-1}- i_1}' \subseteq \ker
N_{i_1}^{b-2}.
\end{gather*} Thus $X_3$ is the component $K_{A_{i_1}}$ of the
Springer fiber $\mathcal{B}_{N_1}$.
\end{proof}

\begin{proof} (of Theorem \ref{hookchar}) A typical fiber $X_2$ of the
map $p_2:K_A \to Fl(F_{i_1})$ consists of flags $F \in K_A$ with fixed
subspaces $F_1, \dots, F_{i_1}$. This fiber $X_2$ is in bijection with
the set of flags in $V/{F_{i_1}}$ that are fixed by $N_{i_1}:
V/F_{i_1} \to V/F_{i_1}$.

 So we have exhibited the component $K_A = X_1$ of the Springer fiber
 $\mathcal{B}_N$ as the
total space of a bundle $p_1: X_1 \to B_1$ with base $B_1 = G_{i_1
-1}(\ker N / \im N^{b-1})$.  The fiber $X_2$ of $p_1$ is the total
space of another bundle $p_2: X_2 \to B_2$ with base
$B_2=Fl(V/F_{i_1})$.

Successive applications of Lemma \ref{smallerspringerlem} prove that
if $X_{2j+1}$ is a component of the Springer fiber of $N_{i_j}:
V/F_{i_j} \to V/F_{i_j}$, then $X_{2j+3}$ is the corresponding
component of the Springer fiber for $N_{i_{j+1}}: V/F_{i_{j+1}} \to
V/F_{i_{j+1}}$.

Finally, since $\im N \subset F_{i_{b-1}}$, we see that the map
$N_{i_{b-1}}: V/F_{i_{b-1}} \to V/F_{i_{b-1}}$ is the zero map. Thus
the (unique) component of the Springer fiber for $N_{i_{b- 1}}$ (which
is $X_{2b-1}$) is the flag manifold $Fl(V/F_{i_{b-1}})$.

\end{proof}

\begin{theorem} Let $N$ be a nilpotent map of hook type, and let $A$
be a standard tableau on the hook shape of $N$, with top row $n,
i_{b-1}, \dots, i_1$. Then the component $K_A$ of the Springer fiber
$\mathcal{B}_N$ has intersection homology Poincar\'e polynomial equal
to
\[[i_1]!\binom{[n-b]}{ [i_1-1]} [i_2
-i_1]!\binom{[(n-i_1)-(b-1)]}{[i_2 - i_1-1]}\dots
[i_{b-1}-i_{b-2}]!\binom{[(n-i_{b-2}-2]}{[i_{b-1}-i_{b-2}
-1]}[n-i_{b-1}]!.\]

This polynomial equals
\[[n-b]![i_1][i_2 - i_1][i_3 -i_2] \dots [i_{b-1} -i_{b-2}][n -
i_{b-1}].\]
\end{theorem}

\begin{proof}

We have already proven that the component in question is an iterated
fiber bundle with $B_{2j}$ a complex flag manifold, with $B_{2j-1}$ a
complex Grassmannian for $1 \leq j \leq b-1$, and with $X_{2b-1} =
B_{2b-1}$ a complex flag manifold. So the bundle $p_{2b-2}: X_{2b-2}
\to B_{2b-2}$ with fiber $X_{2b-1}$ has only even-dimensional homology
in base and fiber. Therefore the Leray-Serre spectral sequence for
$p_{2b-2}$ collapses and the Poincar\'e polynomial of $X_{2b-2}$ is
the product of those for $X_{2b-1}$ and $B_{2b-2}$; in particular,
$X_{2b-1}$ has only even-dimensional homology. Then we do the same
for $p_{2b-3}:X_{2b-3} \to B_{2b-3}$ with fiber $X_{2b- 2}$; since
$B_{2b-3}$ is a complex Grassmannian, the space $X_{2b-3}$ also has
only even-dimensional homology. So $X_{2b-3}, X_{2b-5}, \dots$ have
only even-dimensional homology and their homology Poincar\'e
polynomials are products of Poincar\'e polynomials of flag manifolds
and Grassmannians. This recursion thus unravels to give us the
homology of $X_1$ as stated above. Finally, since the space $X_1$ is
nonsingular, we need only shift the homology Poincar\'e polynomial
until it is invariant under $t \mapsto t^{-1}$, in order to obtain the
intersection homology Poincar\'e polynomial.
\end{proof}

\begin{remark} We have proven that the closed subvariety of
$\mathcal{B}_N$ associated by Vargas' description to the tableau $A$
is irreducible because this subvariety is a bundle of irreducible
varieties; it also contains a Spaltenstein-Vargas set $SV(A)$; hence
it must be exactly the closure $K_A$ of the Spaltenstein-Vargas set
$SV(A)$. This is an alternative proof that Vargas' descriptions indeed
yield the components of the Springer fiber $\mathcal{B}_N$.
\end{remark}

\section{Structure of intersections of two components of hook type}

\begin{theorem} \label{hookintersectionchar} Let $N$ be a nilpotent
map of hook type. Suppose we have two standard tableaux $A$ and $B$ on
the Young shape of $N$, where the standard tableau $A$ has top row $n,
i_{b-1}, \dots, i_1$ and the standard tableau $B$ has top row $n,
i_{b-1}', \dots, i_{1}'$. Then the intersection of the two components
$K_A \cap K_B$ is nonempty iff $\Be_j =\max \{i_j, i_j'\} < \min
\{i_{j+1}, i_{j+1}'\}=\al_{j+1}$, in which case $K_A \cap K_B$ is an
iterated fiber bundle with
\begin{gather*} B_{2j-1} = G_{\Be_j - \Be_{j-1}- 1}(\ker N_{\Be_j}/\im
N^{b({\Be_j})}_{\Be_j})\\ B_{2j} = \{ { F \in
Fl(F_{\Be_j}/F_{\Be_{j-1}})\ |\ \im N_j^{b-j-1} \in F_{\al_j -
\Be_{j-1}} }
\}.\end{gather*}
\end{theorem}

\begin{proof}

We proceed as in the proof of Theorem \ref{hookchar}. Let $\al_j =
\min\{i_j, i_j'\}$ and $\Be_j = \max\{i_j, i_j'\}$. By superimposing
the characterizations of the components $K_A$ and $K_B$, we deduce
that the intersection $K_A \cap K_B$ is given by those flags $F$ in
$Fl(V)$ for which

\begin{gather*}
\im N^{b-1} \subseteq F_{\al_1} \subseteq F_{\Be_1} \subseteq \ker N\\
\im N^{b-2} \subseteq F_{\al_2} \subseteq F_{\Be_2} \subseteq \ker
N^2\\
\im N^{b-3} \subseteq F_{\al_3} \subseteq F_{\Be_3} \subseteq \ker
N^3\\
\dots\\
\im N^{1} \subseteq F_{\al_{b-1}} \subseteq F_{\Be_{b-1}} \subseteq
\ker N^{b-1}.
\end{gather*}

If there exists $j$ for which $\max \{i_j, i_j'\} \geq \min \{i_{j+1},
i_{j+1}'\}$ then the flag $F_{\max \{i_j, i_j'\}}$ would have to
contain $\im N^{b-j -1}$ yet be contained in $\ker N^{j}$, which is
impossible by Lemma \ref{hookimlem}. This proves the emptiness
assertion of the lemma.

Now we exhibit the intersection $K_A \cap K_B = X_1$ as an iterated
bundle.  Define a map $p_1: X_1 \to B_1=G_{\Be_1-1}(\ker N/\im
N^{b-1})$ by $F \mapsto F_{\Be_1}$. The typical fiber $X_2$ of the map
$p_1$ consists of all flags $F$ with a fixed $F_{\Be_1}$ such that
$F_{\al_1}$ contains the $1$-dimensional space $\im N^{b-1}$. Then
there is a map taking $X_2$ to the space $B_2$, which consists of all
flags inside $F_{\Be_1}$ such that $F_{\al_1}$ contains $\im N^{b-1}$,
given by $F_1 \subset F_2 \subset\dots \subset F_n \mapsto F_1 \subset
F_2 \subset\dots
\subset F_{\Be_1}$. We described the structure and homology of the
space $B_2$ in lemma
\ref{IKlem} above.

 Then the fiber of the map $p_2: X_2 \to B_2$ is in bijection with a
 certain space of flags in
$V/F_{\Be_1}$ satisfying (as in the previous theorem) a list of
conditions with respect to the quotient map $N_{\Be_1}: V/F_{\Be_1}
\to V/F_{\Be_1}$. As before, these conditions are exactly the ones
that specify the intersection of two components of the Springer fiber
for $N_{\Be_1}$ whose tableaux $A_{\Be_1}$ and $B_{\Be_1}$ have top
rows $n-\Be_1, i_{b-1} - \Be_1, \dots, i_2 -
\Be_1$ and $n-\Be_1, i_{b-1}' - \Be_1, \dots, i_2' - \Be_1$
respectively.  Thus by descending induction we have our result.
\end{proof}

\begin{corollary} The intersection of the two components in the above
theorem has intersection homology Poincar\'e polynomial equal to
\begin{multline*}
\binom{[n-b]}{[\Be_1-1]}[\Be_1 - \al_1]!\binom{[\Be_1 - 1]}{ [\al_1
-1]}[\al_1]!
\binom{[n-\Be_1 -b+1)]}{[\Be_2 - \Be_1 -1]}[\Be_2 - \al_2]!
\binom{[\Be_2 - \Be_1 -1]}{[ \al_2 - \Be_1 -1]}[\al_2 -\Be_1]!\\ \dots
 \binom{[n-\Be_{j-1} - b +j-1]}{[\Be_j -\Be_{j-1} -1]}[\Be_j -
 \al_j]!\binom{[\Be_j - \Be_{j-1} - 1]}{[\al_j - \Be_{j-1} -1]}[\al_j
 - \Be_{j-1}]! \dots\\
\binom{[ n- \Be_{b-2} -2]}{[\Be_{b-1} - \Be_{b-2} -1]}[\Be_{b-1} -
\al_{b-1}]!\binom{[\Be_{b-1} -
\Be_{b-2} -1]}{[\al_{b-1} -\Be_{b-2} -1]}[\al_{b-1} - \Be_{b-2}]!
[n-\Be_{b-1}]!
\end{multline*} This polynomial equals
\[ [n-b]! [\al_1] [\al_2 -\Be_1][\al_3 - \Be_2]\dots
[\al_{b-1}-\Be_{b-2}] [n-\Be_{b-1}].\]
\end{corollary}

\section{Determination of the topology of components of Springer
fibers for nilpotent maps of two-row type for
$GL_n(\mathbb{C})$}\label{tworowschap}

In this section we study the Springer fibers of nilpotent maps $N$
whose Young shapes have at most two rows.  Thus $N$ has at most two
Jordan blocks.  We will find that the components are iterated bundles with
$\mathbb{CP}^1$ as base spaces, and we will  relate the intersection
homology Poincar\'{e} polynomials of their pairwise intersections to the
inner products of the Kazhdan-Lusztig basis. In doing so, we will
extend some results of Lorist\cite{Lo} on the topology of the
components of two-row shapes with two boxes in the lower row.

Let $N: V\to V$ be a nilpotent map of two-row type.  Recall that a
flag $F$ with subspaces $0 = F_0 \subset F_1 \subset \dots \subset
F_n=V$ is fixed by $N$ if for all $i$, we have $NF_i
\subseteq F_{i-1}$.  Recall that $b$ is defined to be the least
positive integer with $N^b =0$. Similarly let $b(i)$ be the least
positive integer such that $N_i^{b(i)} = 0$.

\begin{definition} Suppose $N_i: V/F_i \to V/F_i$ is a quotient map of
a nilpotent map $N$.  Suppose that the subspace $F_j/F_i$ contains
$\im N_i^k$ but not $\im N_i^{k-1}$. Then we call $\im N^k_i$ the
\emph{lowest image} contained in $F_j/F_i$ and we denote this lowest
image $\im N^k_i$ by $\lowim_i(F_j)$.  (Note that the image of a higher power
of $N_i$ is a smaller subspace of $V/F_i$.) Similarly, we denote the lowest image of
$N$ that is not contained in $F_j$ by $\lowim(F_j)$.
\end{definition}

\begin{lemma} \label{2dlem} Let $N$ be a nilpotent map of two-row
type. Let $F_i$ be a subspace of $V$ such that $F_i \subseteq \im N$
and $NF_i \subset F_i$. Then the quotient space $N^{- 1}(F_i)/F_i$ is
$2$-dimensional.
\end{lemma}
\begin{proof} Since $NF_i \subset F_i$, we see that indeed $F_i
\subset N^{-1}(F_i)$. Then the dimensionality of $N^{-1}(F_i)/F_i$ is
clear from the Jordan form of the map $N$, since $F_i
\subseteq \im N$.
\end{proof}

\begin{lemma} Suppose $N$ is a nilpotent map corresponding to a
two-row Young shape $\tau$. Let $A$ be a standard tableau of shape
$\tau$ with top row $n, i_{b-1}, \dots, i_1$. Then every flag $F$ in
the Spaltenstein-Vargas subset $SV(A)$ defined in \ref{SVdef}
satisfies the following conditions: $F_i \subset N^{-1}(F_{i-1})$ and
$\im N^{b-j} \subseteq F_{i_j}$.
\end{lemma}

\begin{proof} Let the flag $F$ be in the Spaltenstein-Vargas subset
$SV(A)$. Then by construction of $SV(A)$, every subspace clearly
satisfies the first condition.

Now we prove the second condition by inspecting the procedure used to
specify the flag subspaces of $F$. Recall that if $i$ is in the
$c(i)$-th column of $A$, then $F_i/F_{i-1}$ must lie in $\ker N_{i-1}
\cap (\im N_{i-1}^{c(i)-1} - \im N_{i-1}^{c(i)})$.

 Now we show that if $i$ is on the bottom row of $A$ then, for any
 flag $F$ in $SV(A)$,
$\lowim(F_i) = \lowim(F_{i-1})$; in other words, the subspace $F_i$
will never contain a lower image than $F_{i-1}$ contains.  There are
two cases.  First suppose the highest non-zero image $\im
N^{b_{i-1}-1}_{i-1}$ of $N_{i-1}$ is $2$- dimensional. Then
$F_i/F_{i-1}$ cannot exhaust $\im N^{b_{i-1}-1}_{i-1}$. On the other
hand, if the highest image $\im N^{b_{i-1}-1}_{i-1}$ is
$1$-dimensional, then since the number $i$ is not on the top row of
the tableau $A$, the subspace $F_i/F_{i-1}$ must not equal $\im
N^{b_{i-1}- 1}_{i-1}$.

Now note that $F_1, \dots, F_{i_1 -1}$ do not contain wholly any image
of $N$. To stress our line of argument, note that these subspaces
contain the same image of $N$ as $F_0= \{ 0 \}$ does.  Therefore $\im
N^{b-1}_{i_1 -1} \neq 0$.  Then by construction, $F_{i_1}/F_{i_1 -1}$
must contain $\im N_{i_1-1}^{b-1}$. Since $F_{i_1}$ also contains
$F_{i_1 -1}$, we see by Lemma \ref{quotlem} that $F_{i_1}$ must
contain $\im N^{b- 1}$.

Similarly, $F_{i_j}, \dots, F_{i_{j+1} -1}$ all contain $\im N^{b-j}$
and no lower image of $N$, because each of these subspaces is
constructed not to contain the highest image of the previous quotient
map; and, as before, $F_{i_{j+1}}$ must contain $\im N^{b-j-1}$.  Thus
the lemma is proved.
\end{proof}

\begin{remark} Note that the conditions of the lemma are closed and so
are satisfied by the closure of the Spaltenstein-Vargas subset
$SV(A)$, which is the entire component $K_A$ of the Springer fiber
$\mathcal{B}_N$.
\end{remark}

\begin{theorem}
\label{TBthm} Suppose that $N$ is a nilpotent map of two-row type and
that $A$ is a standard tableau on the Young shape of $N$, with top row
$n, i_{b-1}, \dots, i_1$. For any $i$ between $1$ and $n$, denote by
$T(i)$ and $B(i)$ the lengths of the top and bottom rows of the
tableau obtained from $A$ by deleting $1, \dots, i$.  Suppose the flag
$F$ is contained in the Spaltenstein-Vargas subset $SV(A)$. Then the
subspace $F_i$ contains $\im N^{T(i)}$ and is contained in $\im
N^{B(i)}$.
\end{theorem}

\begin{proof} The assertion about $T(i)$ is proven above. As to the
assertion about $B(i)$, we proceed by induction on $i$.  First note that $F_1$ is
contained in $N^{-1}(F_0)$ and therefore must be contained in the
lowest image which has nontrivial intersection with the kernel. By
inspecting the Jordan form of $N$ we see that this image is exactly
$\im N^{B(1)}$.

Now suppose $F_i$ is contained in $\im N^{B(i)}$. Then there are two
possibilities for $F_{i+1}$.  If $i+1$ is on the bottom, then
$F_{i+1}$ is within $N^{-1}(F_i) \subseteq N^{-1}(\im N^{B(i)})$.  Now
by the Jordan form, we see that $B(i)$ cannot be greater than the
power of the lowest image that contains $\ker N$, so $\im N^{B(i)}$
contains $\ker N$. Thus, by lemma \ref{kerlem},
\[N^{-1}(\im N^{B(i)}) = \im N^{B(i)-1} = \im N^{B(i+1)}.\]

Now if $i+1$ is on the top row of the tableau $A$, then by Theorem
\ref{SVthm}, $F_{i+1}/ F_i$ is equal to $(\im N^{T(i+1)} + F_i)/F_i$.
Since $T(i+1) \geq B(i)$, the subspace $F_{i+1}$ is still contained in
$\im N^{B(i)}$.
\end{proof}

\begin{theorem} \label{tworowchar} Let $N$ be a nilpotent map of
two-row type, and let $A$ be a standard tableau on the Young shape of
$A$ with top row $n, i_{b-1}, \dots, i_1$. Then the component $K_A$ of
the Springer fiber $\mathcal{B}_N$ consists of all flags whose
subspaces satisfy the following conditions:
\[ F_i \subset N^{-1}(F_{i-1}) \hbox{ for each $i$, and }
\] if $i$ is on the top row of the tableau $A$ and $i-1$ is on the
bottom row, then
\[ F_i= N^{-1}(F_{i-2});
\] if $i$ and $i-1$ are both in the top row of $A$, then if $F_{i-1} =
N^{-d}(F_{r})$ where $r$ is on the bottom row then
\[ F_i = N^{-d-1}(F_{r-1})
\] and if $F_{i-1} = N^{-d}(\im N^{b-i})$ where $0 \leq i < n-b$ then
\[ F_i= N^{-d}(\im N^{b-i-1}).
\]

The subspaces that are specified as inverse images of
other spaces will be called
\emph{dependent}; note that they are exactly the subspaces whose
indices are on the top row of the tableau $A$. The other subspaces are
called \emph{independent}.
\end{theorem}

\begin{proof} Denote by $K(A)$ the closed subset of flags that satisfy
the conditions of the theorem. Note that $K(A) \subseteq
\mathcal{B}_N$. Let $F$ be a flag in the Spaltenstein-Vargas subset
$SV(A)$. Then we prove by induction on $i$ that each subspace $F_i$ of
$F$ satisfies the above conditions, so that the flag $F$ lies in
$K(A)$.  Then we will show below that the closed subset of the theorem
is in fact irreducible and of the same dimension as the (nonempty)
subset $SV(A)$. Thus $K(A)$ is exactly the closure of the
Spaltenstein-Vargas subset $SV(A)$ and is thus the component $K_A$ of
the Springer fiber $\mathcal{B}_N$.

First we settle the $i=1$ case. Suppose that the number $1$ is in the
bottom row of $A$.  Then for all $F$ in the Spaltenstein-Vargas subset $SV(A)$, it is the case that $F_1$ must be in $\ker N$, 
which is exactly $N^{-1}(F_0)$. If $1$ is in the top row then $F_1$ must equal the highest image $\im
N^{b-1}$; note that $F_0 = 0 = N^{-0}( \im N^b)$ so $F_1 = N^{-0}(\im
N^{b-1})$.  Also note that the highest image must indeed be
$1$-dimensional in order for $1$ to be in the top row. This settles
the $i=1$ case.

Now suppose that $F_1, \dots, F_i$ satisfy the conditions of the
theorem. Now the number $i+1$ is either on the bottom row of the
tableau $A$, or on the top row of $A$.  If $i+1$ is on the bottom row
of $A$, then the Spaltenstein-Vargas procedure requires only that
$F_{i+1}/F_i$ be contained in $\ker N_i$, which by lemma
\ref{invkerlem} equals $N^{-1}(F_i)$.

Now suppose $i+1$ is on the top row of the tableau $A$. Then the
number $i$ is either on the bottom row of $A$ or on the top row. First
we will prove that if $i$ is on the bottom row, then $F_{i+1} =
N^{-1}(F_{i-1})$.

Recall that $F_i \subset N^{-1}(F_{i-1})$, and $F_{i+1}/F_i$ must be
the $1$-dimensional subspace of $V/F_i$ which is the highest nonzero
image of $N_i$. (Since $i$ was on the bottom row, it is clear that the
tableau obtained by deleting $1, \dots, i$ is not rectangular so the
highest image of $N_i$ is not $2$-dimensional.)

So $F_i$ is in $N^{-1}(F_{i-1})=N_{i-1}^{-1}(0+F_{i-1})$ and by
construction $F_i/F_{i-1}$ cannot be all of the highest nonzero image
of $N_{i-1}$. Therefore there must be other vectors $v \in
N^{-1}(F_{i-1})$ such that $v +F_i \neq 0+F_i$, and also $v +F_{i-1}$
is in the highest image of $N_{i-1}$, hence $v$ is in the highest
image of $N$ which is in not in $F_{i-1}$. Then any such $v$ must be
in the highest image of $N$ which is not contained in $F_i$. Since
$N^{-1}(F_{i- 1})/F_i$ is $1$-dimensional, this proves that
$N^{-1}(F_{i-1})/F_i$ must be the highest image of $N_i: V/F_i \to
V/F_i$.

Thus we have proven that if $i$ is in the bottom row and $i+1$ is in
the top row, then $F_{i+1} = N^{-1}(F_{i-1})$.

Now suppose $i+1$ is on the top row, and $i$ is also on the top
row. Then by induction either $F_i = N^{-d}(F_r)$ where $r$ is on the
bottom row, or else $F_i = N^{-d}(\im N^{a})$ (where $n- b < a \leq
b$).  If $r$ is on the bottom row then, by Theorem \ref{TBthm}, the
lowest image that $F_r$ contains is exactly the same as the lowest
image that $F_{r-1}$ contains. Therefore $N^{-1}(F_{r-1})$ contains
one lower image and is of dimension exactly one larger than that of
$F_r$. Thus $N^{-d- 1}(F_{r-1})$ contains exactly one lower image and
is of one larger dimension than $N^{-d}(F_r)$.

Otherwise $F_i = N^{-d}(\im N^a)$ for some $a > n-b$. Since $a > n-b$,
the dimension of $\im N^{a-1}$ is exactly one larger than the
dimension of $\im N^{a}$.  Thus $N^{-d}(\im N^{a-1})$ contains one
lower image and is one dimension larger than $\im N^{a}$. So
$N^{-d}(\im N^{a- 1})/F_{i}$ equals the highest nonzero image $\im
N_{i}^{b(i)}$ and so $F_{i+1} = N^{-d}(\im N^{a- 1})$. Finally, we
prove below that $K(A) = \clo{SV(A)} = K_A$.
\end{proof}

\begin{proposition} The closed subset $K(A)$ is irreducible, and 
$K(A)$ is an iterated bundle of base type $(\mathbb{CP}^1, \dots, \mathbb{CP}^1$, where there are as many terms as there are
number in the bottom row of the tableau $A$.
\end{proposition}
\begin{proof} Suppose that the shape of $N$ has exactly two rows (if
it has only one row, then $K_A$ is a point).  Let $F$ be a flag in the component $K_A$ of
the Springer fiber $\mathcal{B}_N$. Suppose $F_{j_1}$ is the smallest
independent subspace of the flag $F$. Then $F_{j_1-1}$ is some fixed
subspace of $V$ (necessarily in $\im N$), and $F_{j_1}/F_{{j_1}-1}$
can be any point in the fixed space
$\mathbb{P}(N^{-1}(F_{{j_1}-1})/F_{{j_1}-1}) = \mathbb{CP}^1$. Set
$B_1 = N^{-1}(F_{j_1- 1})/F_{j_1-1}$. The map $p_1: K_A \to B_1$ given
by $F \mapsto F_{j_1}/F_{j_1-1}$ is then a fiber bundle (it is clearly
a proper submersion). The typical fiber $X_2$ of the map $p_1$
consists of all flags $F$ in $K_A$ with the subspace $F_i$ fixed, as
well as with all subspaces of $F$ that are dependent on $F_{j_1}$
fixed. Now find the smallest independent subspace $F_{j_2}$ in
$X_2$. Again, all subspaces smaller then $F_{j_2}$ are dependent, and
thus fixed. So we see that $F_{j_2}/F_{j_2 - 1}$ can be any point in
the fixed space $\mathbb{P}(N^{-1}(F_{j_2 -1})/F_{j_2 -1})= B_2$,
which defines the bundle projection $p_2: X_2 \to B_2$, with fiber
$X_3$. We continue until all independent subspaces are exhausted and
the fiber consists of one flag with all subspaces fixed.
\end{proof}

\begin{theorem} \label{bundlethm} Every component $K_A$ of the
Springer fiber for a nilpotent map of two-row type is an iterated
bundle of base type $(\mathbb{CP}^1, \dots,\mathbb{CP}^1$,
where there are as many terms as there are numbers in the bottom row
of the tableau $A$.
\end{theorem}

\begin{proof} The closed subset $K(A)$ is irreducible, contained in
$\mathcal{B}_N$, and clearly has dimension $n-b$, which is the
dimension of the nonempty Spaltenstein-Vargas subset
$SV(A)$. Therefore, by standard algebraic geometry, $K(A)$ must equal
the component $K_A =\clo{SV(A)}$.
\end{proof}

\begin{example} Consider the standard tableau
\[
\begin{matrix} 5 &4&1 \\ 3 & 2 & \end{matrix} \] Then the
corresponding component of $\mathcal{B}_N$ is the set of flags with
$F_i \subset N^{- 1}(F_{i-1})$,  and the conditions $F_1=im N^2$ and $F_4=N^{-1}(F_2)$, which we can write compactly as  
\[ \im N^2 \subset F_2 \subset F_3 \subset N^{-1}(F_2) \subset V. \]
Similarly, the standard tableau
\[ \begin{matrix}5 & 3 & 2 \\ 4 & 1 & \end{matrix}\] corresponds to
the component
\[ F_1 \subset N^{-1}(F_0) \subset N^{-1}(\im N^2) \subset F_4 \subset
V.\]
\end{example}

\begin{example}
Consider the standard tableau
\[
\begin{matrix} 5 &4&3 \\ 2 & 1 & \end{matrix} \] Then the
corresponding component of $\mathcal{B}_N$ is the set of flags with
$F_i \subset N^{- 1}(F_{i-1})$ and
\[
 F_1 \subset F_2 \subset N^{-1}(F_1) \subset N^{-2}(F_0) \subset V. \]

Similarly, the standard tableau
\[ \begin{matrix}5 & 4 & 1 \\ 3&  2& \end{matrix}\] corresponds to
the component
with $F_i \subset N^{- 1}(F_{i-1})$ and
\[ \im N^2 \subset F_2 \subset F_3 \subset N^{-1}(F_2) \subset V. \]

Their intersection is 
the set of flags with $F_i \subset N^{- 1}(F_{i-1})$ and
 \[ \im N^2 < N^{-1}(F_0) < N^{-1}(\im N^2) < N^{-2}(F_0) < V.\] In
particular, the intersection is not empty, in contrast to the
assertions of Wolper\cite{Wo}.

\end{example}

\begin{example} Consider the standard tableau
\[
\begin{matrix} 10 & 9 & 8 & 7 & 4 & 3 \\ 6 & 5 & 2 & 1 & &
\end{matrix}
\] Then the corresponding component of $\mathcal{B}_N$ is the set of
flags with $F_i \subset N^{-1}(F_{i-1})$ and
\begin{multline*} F_1 \subset F_2 \subset N^{-1}(F_1) \subset
N^{-2}(F_0) \subset F_5 \subset F_6 \subset \\ N^{-1}(F_5) \subset
N^{-4}(F_0) \subset N^{-4}(\im N^5) \subset V.
\end{multline*}
\end{example}

\begin{example} We show how to recover Lorist's description of the
structure of components of Springer fibers of $2$-regular nilpotent
maps (that is, those nilpotent maps $N$ whose Young shapes have two
boxes in the second row). There are two types of components $K_A$,
corresponding to whether the numbers on the bottom row of $A$ are
consecutive (yielding a non-trivial bundle) or not consecutive
(yielding a trivial bundle). If the numbers in the bottom rows $A$ are
not consecutive, say $i < j$, then we get
\begin{multline*}
\im N^{b-1}\subset
\im N^{b-2} \subset \im N^{b-i +1} \subset F_i \subset N^{-1}(\im
N^{b-i+1}) \subset \dots\\
\subset N^{-1}(\im N^{b-j+3}) \subset F_j \subset N^{-2}(\im
N^{b-j+3}) \subset N^{-2}(\im N^{b- j+2}) \subset \dots \subset F_n.
\end{multline*}

If the numbers are consecutive, say $i, i+1$, then we get
\[\im N^{b-1}\subset
\im N^{b-2} \subset \im N^{b-i +1} \subset F_i \subset F_{i+1} \subset
N^{-1}(F_i)\subset N^{-2}(\im N^{b-i +1})\subset \dots \subset F_n.
\]\end{example}

Finally, we derive the scholium that the ``dependence" on one subspace
on another can be thought of as symmetric: if $F_i = N^{-d}(F_r)$ in
the above theorem, then indeed $F_i$ also determines $F_r$; so given
either subspace, we can obtain the other. So ``independent" can be
thought of as ``smallest in the chain of dependencies."
\begin{proposition} \label{symdepprop} Suppose $F_i$ is specified as
 $N^{-d}(F_r)$ in Theorem \ref{tworowchar}. Then the map $N^{d}: F_i
 \to F_r$ is surjective. Thus the subspace $F_i$ determines the
 subspace $F_r$.
\end{proposition}
\begin{proof} In general, $N(N^{-1}(W)) =W \cap \im N$. Thus we need
only ensure that if $F_i$ is specified as $N^{-d}(F_r)$ then in fact
$F_r \subseteq \im N^{d}$. First note that the process described in
the proof of Theorem \ref{tworowchar} never takes an inverse image of
a subspace $F_i$ unless $F_i \subset \im N$. Furthermore, if $F_i =
N^{-d}(F_r)$ where $F_r$ is independent and $r>0$, then (for $F$ in
$SV(A)$) the subspace $F_{r-1}$ is contained in one higher image of
$N$ than $F_{r}$. So if $N^{-d+1}(F_r) \subset \im N^{k}$, then also
$N^{-d}(F_{r-1})\subset \im N^{k}$.  This last statement holds in the
entire component $K_A$. Hence if $F_{i} = N^{-d}(F_r)$ then $F_r
\subseteq
\im N^{d}$, so $N^{d}(F_i) = F_r$.
\end{proof}

\section{Relationship with Kazhdan-Lusztig theory}

Let $W$ be a Coxeter group with simple reflections $S$. Denote the
 Chevalley-Bruhat order by $<$.
 Recall \cite{KL}\cite{Hu} that the Kazhdan-Lusztig
 construction yields elements $C_w'$ in the Iwahori-Hecke algebra of
 $W$,  which give
 distinguished bases for certain representations of the Iwahori-Hecke
 algebra, called left cell
 representations, which are associated to certain subsets
 $\mathcal{C}$ of $W$ called left cells. In particular, this
 construction yields a distinguished basis 
for each irreducible representation of the Iwahori-Hecke algebra 
$\mathcal{H}_n$  of the symmetric group $S_n$.

Now recall \cite{DJ}\cite{Mu}
  that every irreducible representation
 $M$ of the Iwahori-Hecke
algebra  $\mathcal{H}_n$  possesses a
unique (up to a scalar)  nondegenerate symmetric bilinear
form $\la\ ,\ \ra $  that is invariant under  $\mathcal{H}_n$, in the sense that for any $v, v' \in M$, we have $\la T_wv, v' \ra =
\la v, T_w^* v' \ra$. (Here, the involution ${}^*:\mathcal{H}_n \to
\mathcal{H}_n$ is defined by sending $T_w
\mapsto T_{w^{-1}}$ and then extending linearly.)

Thus we can consider the inner products of the Kazhdan-Lusztig
basis vectors of an irreducible representation with respect to this
inner product.  We find that they satisfy equations very reminiscent of a 
possible application of the Beilinson-Bernstein-Deligne-Gabber 
Decomposition Theorem. To
 state these relations, first we determine the eigenvectors and eigenvalues of the
  elements $C_s'$ (for $s$ a simple reflection) 
acting by left multiplication on $\mathcal{H}_n$.

\begin{lemma} \label{heckeeiglem} The eigenvectors of the map
$\mathcal{H}_n \to \mathcal{H}_n$ given by left multiplication by
$C_s'$ are given by
\[
\begin{array}{ll} c_w & \text{where $sw < w$, and} \\ C_s'C_w' -
(t+t^{-1})C_{w}' & \text{where $sw>w$.}\end{array}
\]
\end{lemma}

\begin{proof} The first follows immediately from the formula for 
$C_s'C_w'$ (see \cite{Spr}). To see the second, we calculate \[C_s'C_s'C_w'
- (t +t^{-1})C_s'C_w' = 0.  \] Now note that these eigenvectors span
the Iwahori-Hecke algebra $\mathcal{H}_n$.
\end{proof}

\begin{lemma} \label{eiglem} The eigenvectors for $C_s'$ on a left
cell representation $M_{\mathcal{C}}$ are
\[
\begin{array}{ll} C_w' & \text{where $sw < w$, and} \\ C_s'c_w -
(t+t^{-1})c_{w} & \text{where $sw>w$.}\end{array}
\]
\end{lemma}
\begin{proof} Immediate from Lemma \ref{heckeeiglem}.
\end{proof}

Our equations follow from the following 
\begin{lemma}\label{orthlem} Let $V$ be a representation of an algebra
$A$. Suppose $V$ is equipped with an invariant symmetric bilinear form
$\la\, , \, \ra$. Suppose we have an element $a \in A$ with $a = a^*$,
and let $x$ and $y$ be eigenvectors of $a$ with different eigenvalues.
Then $x$ and $y$ are orthogonal; i.e.  $\la x, y\ra = 0$.
\end{lemma}

\begin{proof} Suppose $x$ has eigenvalue $\lambda$ and $y$ has
eigenvalue $\rho$ under $a \in A$.  Then $\la ax, y\ra = \lambda\la x,
y\ra = \la x, ay\ra = \rho\la x, y\ra $ so if $\lambda
\neq \rho$ then $\la x, y\ra =0$.
\end{proof}

\begin{theorem} Let $\mathcal{C}$ be a left cell yielding a left cell
representation $M_{\mathcal{C}}$. Let $\la\ ,\ \ra $ be an invariant
nondegenerate symmetric bilinear form on $M_{\mathcal{C}}$. Let $s$ be
a simple reflection for $W$. Then for each pair $(x, w)$ with $s
\in L(x)$ and $s \not\in L(w)$ (so $s$ descends $x$ but not $w$), we
have an equation between inner products
\[ (t + t^{-1})\la c_x, c_w\ra = \la c_x, C_s' c_w\ra .\]
\end{theorem}
\begin{proof} Given a simple reflection $s$, we have a large supply of
eigenvectors for $C_s'$ in $M_{\mathcal{C}}$ given by Lemma
\ref{eiglem}, namely $c_x$ where $sx < x$ and $(t + t^{-1}) c_w -
C_s'c_w$ where $sw > w$. Then for each pair $(x, w)$, we have two
eigenvectors of $C_s'$ with different eigenvalues, so the 
eigenvectors are orthogonal. Therefore we have the equations as
claimed.
\end{proof}

In fact, these equations are equivalent to
$\mathcal{H}_n$-invariance: 

\begin{proposition} \label{invarbilinprop} If a bilinear form on a
left cell representation satisfies the equations
\[(t + t^{-1})\la c_x, c_w\ra = \la c_x, C_s' c_w\ra \] for each pair
$c_x$, $c_w$ where $sw > w$ and $sx <x$, then the inner product is
invariant under the action of $C_s'$.
\end{proposition}
\begin{proof} The inner product of any two vectors can be expressed as
a linear combination of inner products of basis vectors, so we are
reduced to proving $\la C_s'c_x, c_w\ra = \la c_x, C_s' c_w\ra $ for
each pair $(x, w)$.

There are three cases. If both $c_x$ and $c_w$ are descended by $s$,
then clearly \[\la (t + t^{- 1})c_x, c_w\ra = \la c_x, (t +
t^{-1})c_w\ra .\]

If $c_x$ is descended by $s$ but $c_w$ is not, then by the equations
above
\[ \la C_s'c_x, c_w\ra = \la (t + t^{-1})c_x, c_w\ra = \la c_x, C_s'
c_w\ra .\]

Finally suppose neither basis vector in $\la c_x, c_w\ra $ is
descended by $C_s'$. Then, recalling that $C_s'c_x = \sum_{\substack{y
\simeq_L x \\ sy < y}} \mu(x, y)c_y$, we have
\[ (t+t^{-1})\la C_s'c_x, c_w\ra = \la C_s' C_s'c_x, c_w\ra = \la
C_s'c_x, C_s'c_w\ra ,
\] since each term in $C_s' c_x$ is descended by $s$.  Then by a
symmetric argument
\[ \la C_s'c_x, C_s'c_w\ra = \la c_x, C_s'C_s'c_w\ra = (t+t^{-1})\la
 c_x, C_s'c_w\ra .\] This completes the proof.
\end{proof}

 The form of these equations is very similar to the
conclusion of the Decomposition Theorem
\cite{BBD}\cite{CG}. This analogy suggests an interpretation in terms of 
Poincar\'e polynomials: $\la c_x, c_w\ra $ as the intersection homology Poincare\'e polynomial of the intersection of two
spaces $K_A$ and $K_B$, $(t + t^{-1})$ as a $\mathbb{CP}^1$, the left
side as a $\mathbb{CP}^1$ bundle over $K_A \cap K_B$, and the right
side as terms from the Decomposition Theorem, applied to some map from
the total space of the bundle to some space.

Thus, given a $W$-graph for a left cell representation, we can write
down a set of equations that the Gram matrix entries of an
  $\mathcal{H}_n$-invariant inner product must
satisfy, and which determine them up to a common scalar.
If we can prove that the Poincar\'{e} polynomials of some collection of 
spaces satisfy those equations, and we can calculate the Poincar\'{e}
 polynomial of one of the spaces, then we can calculate the Poincar\'{e}
polynomials of all the spaces in the collection. We will now do this for 
Springer fibers of hook and two-row type, where the $W$-graphs are known 
explicitly.

The $W$-graphs for left cells with hook shapes are very easily
determined because no two standard tableaux have the same descent
set. (See Garsia-McLarnan\cite{GMcL}, Fact 14, or Curtis-Iwahori-
Kilmoyer\cite{CIK}.) Since all left cell
representations of $\mathcal{H}_n$ for a given Young shape are
isomorphic(\cite{KL}, Theorem 1.4), we can label the Kazhdan-Lusztig basis vectors by their
tableaux.
\begin{definition} A standard tableau of hook shape $B$ is
\emph{adjacent} to the standard tableau $A$ via $k$ if $B$ can be
obtained from $A$ by exchanging $k$ with either $k+1$ or $k-1$ in the
tableau $A$.  Note that exactly one of $A$ and $B$ will have $k$ as a
descent.
\end{definition}

Note that for each standard tableau $A$ and number $k$, there are at
most two standard tableaux adjacent to $A$. We denote by $(k-1\ k)A$
the tableau obtained by switching $k$ and $k-1$ in $A$; similarly,
$(k\ k+1)A$ is the tableau obtained from $A$ by switching $k$ and
$k+1$.

Therefore, in the case of hook shapes, we can explicitly exhibit the
set of equations that the inner products of the Kazhdan-Lusztig basis
vectors must satisfy. Using the above description of the Poincar\'e
polynomials of the intersection homology of the components of the
Springer fibers and their intersections, we will show 
 that the intersection homology Poincar\'e
polynomials of the intersections of the components satisfy the same
equations. Since the equations determine the inner product up to a
scalar by Proposition
\ref{invarbilinprop}, we will be able to show that the inner product
matrix of the Kazhdan-Lusztig basis vectors computes the intersection
homology of the Springer fibers for hook shapes.

\begin{theorem}\label{hookthm} Suppose $A$ is a hook shape tableau with top row $n,
i_{b-1}, \dots, k, \dots, i_1$, so that the number $k$ is not a
descent of the tableau $A$. Then the Kazhdan-Lusztig basis vector
$c_A$ transforms as

\[ T_{(k \ k+1)} c_A = -c_A +t c_{(k\ k-1)A} +t c_{(k\ k+1)A}.
\]If $(k\ k+1)A$ or $(k-1\ k)A$ is not standard, then omit the
corresponding term in the formula above.

If $B$ is a standard tableau that does not have $k$ in the top row
then
\[ T_{(k\ k+1)}c_{B} = t^2 c_B.
\]
\end{theorem}

\begin{proof}

Given two different Young tableaux $A$ and $A'$ of hook type, there
exists a simple reflection that descends $A$ but not $A'$, and a
different simple reflection that descends $A'$ but not $A$, because
their first columns are distinct.  Given a left cell $\mathcal{C}$ of
hook type, there is a $W$-graph for the left cell representation,
indexed by elements of $\mathcal{C}$. 
By (\cite{Hu}, Proposition 7.15), the only edges in the $W$-graph
 are those that
connect elements of the form $x$ and $sx$, with $s$ a simple
reflection.  So in particular the only possible $W$-graph neighbors to
$A$ with $k$ as a descent are $(k-1\ k)A$ and $(k\ k+1)A$. Depending
on which of $k-1, k+1$ can be interchanged with $k$ in the tableau
$A$, we arrive at the above possibilities.
\end{proof}

\begin{lemma} Given the above tableaux $A$, $B$ the following are
eigenvectors for the action of $C_{(k \ k+1)}'$:
\begin{gather*} (t+t^{-1})c_A - c_{(k\ k-1)A} - c_{(k\ k+1)A} \text{
with eigenvalue $0$}\\ c_B \hbox{ with eigenvalue $(t + t^{-1})$}
\end{gather*} If $(k\ k+1)A$ or $(k-1\ k)A$ is not standard, then omit
the corresponding term in the formula above.
\end{lemma}

\begin{proof} This is immediate from Lemma \ref{eiglem} and the
multiplication formula above.
\end{proof}

\begin{theorem} Suppose we are given a nilpotent map $N$ on $V$ with a
Young shape $\tau$ of hook type. Let $TOP$ be the standard tableau on
the shape $\tau$ with top row $n, b-1, b-2,\dots, 1$.  Normalize the
inner products of the Kazhdan-Lusztig basis vectors so that the norm
$\la c_{TOP}, c_{TOP}\ra$ has the intersection homology Poincar\'e
polynomial of the Springer fiber component $K_{TOP}$. Then the inner
product $\la c_{A}, c_{B}\ra$ is equal to the intersection homology
Poincar\'e polynomial of the intersection $K_A \cap K_B$.
\end{theorem}

To accomplish this, we define maps between certain spaces, from which
the  Decomposition Theorem asserts that these intersection homology
 Poincar\'{e} polynomials satisfy the equations of the Kazhdan-Lusztig
inner products. We need some geometric preliminaries.

\begin{definition} A subvariety $X$ of the flag manifold $Fl(V)$ is a
\emph{union of lines of type $k$} if whenever $X$ contains a flag
$F_1\subset\dots \subset F_{k-1} \subset F_{k} \subset F_{k+1}
\subset \dots \subset F_n$, then $X$ contains all flags of the form
$F_1\subset\dots \subset F_{k- 1}\subset F_k' \subset F_{k+1} \subset
\dots \subset F_n$ where $F_k'$ is between the given subspaces
$F_{k-1}$ and $F_{k+1}$. We will, by analogy with Weyl groups, also
say that $k$ is a
\emph{descent} of $X$.  Some authors say that $X$ is
\emph{$k$-vertical}.
\end{definition}
\begin{definition} Let $X$ be a subvariety of the flag manifold of
$GL_n(\mathbb{C})$ and $1\leq k \leq n-1$. We will denote by
$\mathbb{CP}^1\star X$ the variety of pairs
\[
\{ (F_k', F) \text{ where $F \in X$ and $F_k'$ lies between $F_{k-1}$
and $F_{k+1}$}
\}.
\]
\end{definition} This variety admits a map $\phi: \mathbb{CP}^1\star X
\to Fl(V)$ given by mapping $(F_k', F)$ to the flag $F' = F_1 \subset
\dots \subset F_{k-1} \subset F_k' \subset F_{k+1} \subset \dots
\subset F_n$.
\begin{definition} The image $\phi(\mathbb{CP}^1\star X) \subseteq
Fl(V)$ is called the
\emph{$k$-saturation} of $X$; this image is denoted $S_k(X)$.
\end{definition} So the $k$-saturation $S_k(X)$ is obtained roughly by
taking all flags in $X$ and allowing $F_k$ to vary freely within them.
Note that $\phi^{-1}(F)$ is a $\mathbb{CP}^1$ exactly when $F$ is
contained in a line of type $k$ in $X$.

\begin{lemma} The space $\mathbb{CP}^1\star X$ is a locally trivial
fiber bundle over $X$ with fiber $\mathbb{CP}^1$ via the obvious map
$(F_k', F)\mapsto F$.
\end{lemma}
\begin{proof} This is clear.
\end{proof}

\begin{remark} We can also rephrase this using a minimal parabolic
subgroup $P_k$ corresponding to the simple reflection $(k\ k+1)$. Then
the map $\pi: G/B \to G/P_k$ forgets about $F_k$, so the preimage
$\pi^{-1}(\pi(X))$ is equal to $\phi(\mathbb{CP}^1\star X)$.
\end{remark}

Suppose we have two standard tableaux $A$ with top row $n, i_{b-1},
\dots,k=i_j,\dots, i_1$ and $B$ with top row $n, i_{b-1}',
\dots,\hat{k}, \dots, i_1'$ (that is, $k$ is in the $j$th position from the 
right in $A$, but is not in the top row in $B$).  Then $k$ is a descent of the component
$K_B$ but not of the component $K_A$. Let the position of $k$ in the
top row of $A$ be $i_j$. Then we have the adjacent components $K_{(k\
k-1)A}$ and $K_{(k\ k+1)A}$ of the Springer fiber $\mathcal{B}_N$ that
have $k$ as a descent. Then we shall show that the intersection
homology Poincar\'e polynomial of the $k$-saturation of the
intersection $K_B \cap K_A$ equals the sum of the intersection
homology Poincar\'e polynomials of $K_B \cap K_{(k\ k-1)A}$ and $K_B
\cap K_{(k\ k+1)A}$, and we shall show how this equality can be
interpreted in terms of the Decomposition Theorem.

\begin{theorem} Let $N$ be a nilpotent map of hook type. Let $A$ and
$B$ be two standard tableau on the Young shape of $N$ such that $k$ is
a descent of $B$ but not of $A$. Then the intersection $K_A\cap K_B$
of the two components $K_A$ and $K_B$ of the Springer fiber
$\mathcal{B}_N$ is not $k$-saturated. Then the map $\phi$ from
$\mathbb{CP}^1 \star (K_B\cap K_A)$ to the $k$- saturation $S_k(K_B
\cap K_A)$ yields an equation of intersection homology Poincar\'e
polynomials
\[ IP(\mathbb{CP}^1 \star (K_B\cap K_A)) = IP(K_{(k\ k+1)A}\cap K_B )+
IP(K_{(k-1\ k)A}\cap K_B).
\] Remember that the intersection homology Poincar\'e polynomial $IP$
is normalized so that the sum is centered around the degree $0$ term.

If either $(k-1 \ k )A$ or $(k \ k+1)A$ is not standard, then omit the
corresponding term in the above. If both are not standard then the
intersection $K_B \cap K_A$ is empty.

\end{theorem}

\begin{proof} Suppose that we have the intersection $K_A \cap K_B$ of
 two components $K_A$ and $K_B$ where $k$ is a descent of $B$ but not
 of $A$. Let $k =i_j$ in $A$. First, if neither $(k-1 \ k )A$ or $(k
\ k+1)A$ is standard, then we can check by Theorem
\ref{hookintersectionchar} that the intersection $K_A
\cap K_B$ is empty.  In fact, suppose $K_A \cap K_B$ is nonempty. Then
if $A$ has top row $i_{j+1} = k+1$, $i_j =k$ and $i_{j-1} = k-1$, then
$\Be_{j-1} = k-1$ and $\al_{j+1} = k+1$.  Then in the tableau $B$, the
entry $i_j'$ must satisfy $k-1 < i_j' < k+1$, but $i_j' \neq k$
because $k$ is a descent of $B$, which is a contradiction.

Suppose first that $(k \ k+1)A$ is standard. Then $i_j = k$, $i_{j-1}
< k-1$ and $i_{j+1} > k+1$.  Now the term $i_j'$ in $B$ must satisfy
either $i_j' > k$ or $i_j' < k$. Suppose now that $i_j' > k$ (which is
only possible when $(k \ k+1)A$ is standard), so that $k$ = $\al_j$
and $i_j' =
\Be_j$. Then for any $F \in K_A \cap K_B$, the subspaces of $F$
satisfy

\[
\im N^{b-j} \subset F_k \subset F_{k+1} \subset \dots \subset F_{i_j'}
\subset \ker N^{j}
\]

Note that $K_A\cap K_B$ has a (Zariski) open subset of flags with $\im
N^{b-j} \subset F_{k}$ but $\im N^{b-j} \not\subset F_{k-1}$. In such
flags, $F_k$ is determined by $F_{k-1}$ (and $\im N^{b-j}$), so the
intersection $K_A \cap K_B$ is not a union of lines of type $k$.

So the $k$-saturation of the intersection $K_A \cap K_B$ consists of
all flags $\cdots F_{k-1}
\subset F_k' \subset F_{k+1}\cdots$, where $F_k'$ is any subspace
between $F_{k-1}$ and $F_{k+1}$.  In particular, $F_k'$ no longer need
contain $\im N^{b-j}$.  However, the subspace $F_{k+1}$ must still
contain $\im N^{b-j}$ in all of the resulting flags.  Therefore the
$k$-saturation of $K_B \cap K_A$ is all flags with
\[
\im N^{b-j} \subseteq F_{k+1} \subset F_{i_j} \subseteq \ker N^{j},
\] and all the other conditions unaffected.  Therefore the
$k$-saturation $S_k(K_A \cap K_B)$ is indeed $K_{(k\ k+1)A} \cap K_B$.

If $(k-1\ k)A$ is also standard, then there will be a nonempty subset
of $K_B \cap K_A$ consisting of those flags in the intersection for
which
\[
\im N^{b-j} \subset F_{k-1} \subset F_{i_j} \subset \ker N^{j}
\] which clearly corresponds to the intersection $K_B \cap K_{(k-1\
k)A}$; this subset is a union of lines of type $k$. The subvariety
$K_B \cap K_{(k-1\ k)A}$ is of codimension $2$ in $K_B \cap K_{(k\
k+1)A}$.

Therefore we have a map $\phi: \mathbb{CP}^1 \star (K_B\cap K_A) \to
K_{(k\ k+1)A}\cap K_B$. This map is generically $1-1$. The map $\phi$
is a semismall resolution, because the subvariety where $\phi$ has a
$\mathbb{CP}^1$ fiber is exactly $K_{(k-1\ k)A}\cap K_B$, which is of
codimension $2$ in the image space $K_{(k\ k+1)A}\cap K_B$, and the
domain is nonsingular.

 Therefore when we invoke the Decomposition Theorem for semismall maps
 \cite{BM}\cite{CG}, we find that the
 intersection homology
Poincar\'e polynomial $IP(\mathbb{CP}^1 \star (K_B\cap K_A))$ equals
the intersection homology Poincar\'e polynomial $IP(K_{(k\ k+1)A}\cap
K_B)$ of the range, plus the intersection homology $IP(K_{(k-1\
k)A}\cap K_B)$ of the smaller intersection. (Remember that the
intersection homology Poincar\'e polynomial $IP$ is normalized so that
all of these sums will be centered around $0$.)

If $(k\ k+1)A$ is standard but $(k-1 \ k)$ is not, then $K_{(k-1\
k)A}\cap K_B$ is empty, so its term is omitted. Finally, if $i_j' < k$
so $k=\Be_j$ (which is only possible if $(k-1\ k)A$ is standard), then
the roles of $K_{(k-1\ k)A}$ and $K_{(k\ k+1)A}$ will be reversed in
the above argument.
\end{proof}

\begin{example} We have an equality of $IP$'s of the following spaces (where
we denote the space $K_A$ by its tableau $A$):
\begin{multline*}
\mathbb{CP}^1\star\left(
\begin{matrix}7&6 & 3 & 2 \\ 5&&&\\ 4&&&\\ 1&&& \end{matrix}
\quad \bigcap
\quad
\begin{matrix}7&5 & 3 & 1 \\ 6&&&\\ 4&&&\\ 2&&& \end{matrix}
\right)
\\ =
\left(
\begin{matrix}7&6 & 3 & 2 \\ 5&&&\\ 4&&&\\ 1&&& \end{matrix}\quad
\bigcap \quad
\begin{matrix}7& 6& 3 & 1 \\ 5&&&\\ 4&&&\\ 2&&& \end{matrix}
\right)
\qquad {+} \qquad
\left(
\begin{matrix}7&6 & 3 & 2 \\ 5&&&\\ 4&&&\\ 1&&& \end{matrix}\quad
\bigcap \quad
\begin{matrix}7&4 & 3 & 1 \\ 6&&&\\ 5&&&\\ 2&&& \end{matrix}
\right)
\end{multline*}

For the first intersection, we have
\[\al_1 =1, \Be_1 = 2, \al_2=3, \Be_2 =3, \al_3 = 5, \Be_3 = 6\] so
the computation of the polynomials is
\begin{multline*}[2] \times [1][7-4]![3-2][5-3][7-6] \\ =
[1][7-4]![3-2][6-3][7-6] + [1][7-4]![3-2][4-3][7-6].\end{multline*}
\end{example}

\section{Structure of intersections of components of two-row type and
the relationship with Kazhdan-Lusztig theory} Now we compute the
intersection homology Poincar\`e polynomials of pairwise intersections
of components and then relate them to Kazhdan-Lusztig theory. In doing
so, we also correct a result of Wolper\cite{Wo}.

Let $N$ be a nilpotent map with a two-row Young shape $\tau$.  In this
section we prove that the intersection homology Poincar\'e polynomials
of the intersections of the components of the Springer fiber
$\mathcal{B}_N$ coincide with the (appropriately normalized) inner
products of the Kazhdan-Lusztig basis vectors of the left cell
representation of $\mathcal{H}_n$ associated to the Young shape
$\tau$.

The Kazhdan-Lusztig inner product matrix has been studied for left
cell representations of two- row type because they yield the
representations of the so-called Temperley-Lieb algebra\cite{We}.  We
shall understand how the combinatorics of the Temperley-Lieb algebra
representations encodes the structure of the intersections of
components of the Springer fiber.

First, let us review some notions from Temperley-Lieb theory\cite{We}.

\begin{definition} Suppose we have the numbers $1$ to $n$ on a
horizontal line, increasing to the right. Then an $(n, p)$-\emph{cup
diagram} consists of $p$ \emph{cups} on these numbers, where each cup
connects two numbers, no two cups intersect each other, and no number
is underneath a cup and yet not connected to any cup. The entire cup
diagram must lie in one half- plane.
\end{definition}

\begin{lemma} Suppose $\tau$ is a Young shape with $n$ boxes and $n-b$
boxes in the second row. Then there is a bijection between standard
Young tableaux on $\tau$ and $(n, n-b)$-cup diagrams. This bijection
is denoted $A \to CupD(A)$.
\end{lemma}

\begin{proof} Let $A$ be a standard tableau of shape $\tau$. We
construct a cup diagram as follows. Begin at the number $1$ on the
horizontal line. Proceed from left to right, starting a cup at $i$ if
the number $i$ is on the bottom row of $A$, and ending a cup if $i$ is
on the top row by matching the number $i$ with the closest started cup
that can be matched with $i$. All unpaired ends of cups are then left
blank (these are called \emph{orphaned} numbers). It is easily seen
that this procedure produces a bijection between cup diagrams and
two-row tableaux.
\end{proof}

\begin{example} The $(7,3)$ cup diagram
\vskip12pt
\centerline{\epsfbox{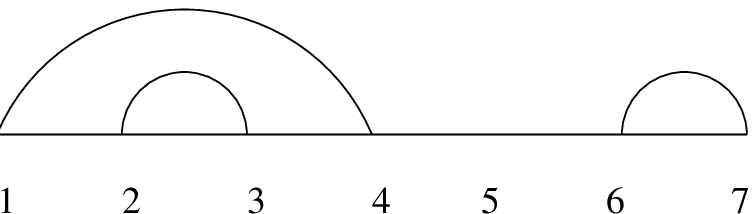}}
\vskip12pt corresponds to the standard Young tableau

\[ \begin{matrix} 7&5&4&3 \\ 6&2&1& \end{matrix}\]
\end{example}
\begin{lemma} A two-row standard tableau $A$ has a descent at $i$ if
and only if the associated cup diagram $CupD(A)$ has a cup connecting
$i$ and $i+1$. Such a cup will be called a \emph{minimal} cup.
\end{lemma}
\begin{proof} The tableau $A$ has a descent exactly when $i$ is on the
bottom row and $i+1$ is on the top row.  This is the case exactly when
there is cup connecting the numbers $i$ and $i+1$ in the cup diagram
$CupD(A)$.
\end{proof}

We shall exhibit a correspondence between the
 cup diagram $CupD(A)$ and the dependencies among the subspaces of the
 flags in the component $K_A$. First, let us extend the cup diagram
 $CupD(A)$ by adding the numbers $0$, $-1$, $\dots$, $-(n-b)$ to the
 left of the numbers $1$, $\dots$, $n$. Now match each orphaned number
 in the cup diagram, working from left to right, to the closest
 possible negative number. This creates an extended cup diagram that
 we label $ECupD(A)$. Note that there are no additional choices here
 so these extended cup diagrams are still in bijection with the
 standard tableaux. For a point $i$ which is an endpoint of a cup,
 denote by $\sigma(i)$ the other endpoint of the cup.

\begin{definition} Given a cup diagram $CupD(A)$ and an index $i$ at
 which a cup begins or ends, we denote by $Cup_A(i)$ the cup that
 begins or ends at $i$.
\end{definition}

Recall that a subspace $F_i$ in a component $K_A$ is called
\emph{dependent} if $F$ is specified by Theorem \ref{tworowchar} as
the inverse image of some smaller subspace (equivalently, the number
$i$ is on the top row of $A$), and is said to
\emph{depend} on that smaller space.  Otherwise the subspace is called
\emph{independent}.

\begin{remark} \label{negconv} For $-(n-b) \leq i < 0$, we
interpret $F_i$ to mean $\im N^{b-i}$.
\end{remark}

\begin{definition} A cup $Cup_1$ lies \emph{directly beneath} a cup
$Cup_2$ if $Cup_1$ is beneath $Cup_2$ and there are no other cups that
lie both above $Cup_1$ and below $Cup_2$.
\end{definition}

\begin{theorem} Consider a nilpotent map $N$ and a two-row standard
tableau $A$ on the Young shape of $N$. Then the extended cup diagram
$ECupD(A)$ encodes the dependencies among the subspaces of the flags
in $K_A$ as follows.  If a cup begins at $i$, then $F_i$ is
independent. If a cup ends at $i$ then $F_i$ is an inverse image of
$F_{\sigma(i) - 1}$ (interpreted using convention
\ref{negconv}). \end{theorem}

\begin{proof} Let $i>0$. First, note that $F_i$ is independent iff the
number $i$ is on the bottom row of the tableau $A$. This is true iff
$i$ starts a cup in $ECupD(A)$. The subspace $F_i$ is dependent on a
smaller subspace iff $i$ ends a cup in $ECupD(A)$.

Now suppose $F_i$ is dependent. We apply the characterization of
Theorem \ref{tworowchar}.  Let us proceed by induction on the length $|\sigma(i) -i
|$ of the cup $Cup_A(i)$.  If $i-1$ is independent, then $F_{i}$ is
equal to $N^{-1}(F_{i-2})$. On the other hand, $i-1$ starts a cup and
$i$ ends the cup, so there must be a minimal cup connecting $i-1$ and
$i$, and so $i-2 = \sigma(i) -1$. This proves the $|\sigma(i) - i| =
1$ case.

Now suppose that $i$ is dependent and $i-1$ is also dependent. Then,
since cups cannot cross, we see that $\sigma(i) < \sigma(i-1)$, so
that the cup $Cup_A(i-1)$ is shorter than the cup $Cup_A(i)$. Now note
that all numbers under the cup $Cup_A(i)$ must either begin or end a
cup.  So if $\sigma(i-1) \neq \sigma(i) +1$, then there must exist a
sequence of adjacent cups directly beneath $Cup_A(i)$ beginning at
$\sigma(i)+1$ and ending at $\sigma(i-1)-1$.  Then by induction, the
space $F_{i-1}$ depends on the independent subspace $F_{\sigma(i)}$,
so by Theorem \ref{tworowchar} the subspace $F_i$ depends on the
subspace $F_{\sigma(i) - 1}$.
\end{proof}

\begin{proposition} \label{cupdependprop} Suppose $i>0$ begins a cup
in $CupD(A)$.  If a subspace $F_j$ depends on the (independent)
subspace $F_i$, then $i < j < \sigma(i)$; that is, $j$ lies strictly
under the cup that begins at $i$. In fact, if $j$ is the end of a cup
that lies directly beneath $Cup_A(i)$ then $F_j$ depends in $F_i$.
\end{proposition}
\begin{proof} Note first that $j >i$ since $F_j$ is to depend on the
independent subspace $F_i$. Recall that since $F_j$ is dependent, $j$
ends a cup, and $F_j$ depends on the subspace $F_{\sigma(j) -1}$.
Either that subspace is independent, or $\sigma(j) - 1$ ends another
cup so $F_{\sigma(j) -1}$ depends on $\sigma(\sigma(j) - 1) -1$, and
so forth.

If $j > i$ and $j$ does not lie strictly under the cup starting at $i$
then either $\sigma(j) < i$ or $\sigma(j) >
\sigma(i)$ so $\sigma(j) - 1$ cannot lie strictly under the cup
either. If $\sigma(j) - 1 =
\sigma(i)$ then $F_j$ depends on $F_{i-1}$ and thus not on $F_i$.  So,
$F_{\sigma(j) - 1}$ cannot depend on $F_i$ unless $i < j < \sigma(i)$.

As to the last assertion, note that if $j$ ends a cup lying directly
beneath $Cup_A(i)$, then $F_j$ must depend on $F_{\sigma(j) -1}$. Then
$\sigma(j) - 1$ must also end a cup lying directly beneath $Cup_A(i)$
unless $\sigma(j) - 1 = i$. This completes the proof.
\end{proof}

Now we demonstrate that the intersections of components of the
Springer fiber $\mathcal{B}_N$ satisfy the equations for the
Kazhdan-Lusztig inner products. Suppose $A$ is a standard tableau on
the shape of $N$. Suppose $i$ is not a descent of $A$ (so it is not
the case that $i$ is on the bottom row and $i+1$ is on top). The
assertion that $i$ is not a descent is equivalent to the assertion
that there is not a cup joining $i$ and $i+1$. Then, we can
manufacture a cup diagram having $i$ as a descent.

\begin{definition} Suppose $CupD(A)$ is a cup diagram that does not
have a minimal cup connecting $i$ and $i+1$.  Suppose $\sigma(i) \neq
i$ and $\sigma(i+1) \neq i+1$.  Then the cup diagram $CupD(u_i A)$ is
defined by deleting the cups with endpoints at $i$ and $i+1$, then
connecting $i$ and $i+1$ with a minimal cup and connecting $\sigma(i)$
to $\sigma(i+1)$ with another cup.  If exactly one of $\sigma(i) = i$
or $\sigma(i+1) = i+1$, then we only insert the cup between $i$ and
$i+1$.  If both $\sigma(i) = i$ and $\sigma(i+1) = i+1$ then $CupD(u_i
A)$ does not exist. Note that this definition also defines a standard
tableau $u_i A$.
\end{definition}

In Westbury\cite{We} it is proven that the tableau $u_i A$ gives the
unique $W$-graph neighbor to $A$ that has $i$ as a descent; if this
tableau $u_i A$ does not exist, then there are no neighbors to $A$ in
the $W$-graph with $i$ as a descent. We now show that if we
$i$-saturate the intersection $K_A \cap K_B$ (where $K_B$ is descended
by $i$) then we get the intersection $K_{u_i A} \cap K_B$.

\begin{theorem}\label{tworowthm} Let $N$ be a nilpotent map of two-row type. Consider
two standard tableaux $A$ and $B$ on the Young shape of $N$ such that
$i$ descends $B$ but not $A$. Suppose $u_i A$ is the unique $W$-graph
neighbor to $A$ that has $i$ as a descent. Then the intersection
homology Poincar\'e polynomials of the intersections satisfy the
following equality:
\[ (t +t^{-1})IP(K_A \cap K_B) = IP(K_{u_i A} \cap K_B).\] If there is
no such neighbor $u_i A$ then the intersection $K_A \cap K_B$ is
empty.
\end{theorem}

\begin{proof} We shall show that the $i$-saturation
$\mathbb{CP}^1\star (K_A \cap K_B)$ of the intersection $K_A
\cap K_B$ has $F_i$ independent, but all other dependencies among the
other subspaces are the same as in $K_A \cap K_B$. This will prove the
theorem. We use the fact that if $F_j$ depends on $F_i$, then we can
also determine $F_i$ from knowledge of $F_j$ (Proposition
\ref{symdepprop}).

There are several cases; first note that in all cases, $F_{i+1} =
N^{-1}(F_{i-1})$ because this dependency holds in $K_B$.

\begin{enumerate}
\item Suppose $\sigma(i+1) < \sigma(i) < i < i+1$ in $A$.  First, for
a flag $F$ in the component $K_A$, we see that $F_{i+1}$ depends on
$F_{\sigma(i+1) -1}$.  Because $F_{i+1} = N^{-1}(F_{i-1})$ in $K_B$,
we see that $F_{i+1}$ must depend on $F_{i-1}$.  Now in the component
$K_A$, the number ${i-1}$ is either equal to $\sigma(i)$, or $i-1$
lies at the end of a cup lying directly under $Cup_A(i)$. So the
subspace $F_{i-1}$ depends on the subspace for the start of the cup
$Cup_A(i)$, namely $F_{\sigma(i)}$. Thus $F_{i+1}$ depends on
$F_{\sigma(i)}$. Since $F_{i+1}$ depends on $F_{\sigma(i+1)-1}$ in
$K_A$, we see that $F_{\sigma(i)}$ must depend on $F_{\sigma(i+1)-1}$
in the intersection $K_A \cap K_B$.

Then in the transformed tableau $u_i A$, we now have a cup connecting
${\sigma(i+1)}$ to ${\sigma(i)}$ and one connecting $i$ to $i+1$.
This means that for any flag $F$ in the component $K_{u_i A}$,
$F_{\sigma(i)}$ is dependent on $F_{\sigma(i+1) -1}$. Thus the
intersection $K_{u_i A} \cap K_B$ will have all the same dependencies
between subspaces as $K_A
\cap K_B$ does, except that $K_{u_i A} \cap K_B$ will be
$i$-saturated.

\item Suppose $\sigma(i)< i < i+1 < \sigma(i+1)$. Then in $K_A$, the
subspace $F_i$ depends on $F_{\sigma(i)-1}$, and $F_{\sigma(i+1)}$
depends on $F_i$ and thus on $F_{\sigma(i)-1}$.

Now for any flag $F$ in $K_{u_i A}$, the subspace $F_{\sigma(i+1)}$
depends on $F_{\sigma(i) -1}$ as well. No dependency is imposed on
$F_{i+1}$ that was not present in $K_A \cap K_B$.

\item Suppose $i < i+1 < \sigma(i+1) < \sigma(i)$. Then for any flag
$F$ in the component $K_A$, the subspace $F_{\sigma(i)}$ depends on
$F_{i-1}$, and the subspace $F_{\sigma(i+1)}$ depends on $F_i$. Also
${\sigma(i+1) - 1}$ is either $i+1$ or is the end of a cup directly
under $Cup_A(i+1)$. So, by Proposition \ref{cupdependprop}, the
subspace $F_{\sigma(i+1) - 1}$ depends on $F_{i+1}$. Thus in the
intersection $K_A \cap K_B$, the subspace $F_{\sigma(i+1) - 1}$
depends on $F_{i-1}$.

Then in $K_{u_i A}$, we have that $F_{\sigma(i)}$ depends on
$F_{\sigma(i+1) -1}$. By the same chain of dependencies,
$F_{\sigma(i+1) -1}$ depends on $F_{i+1}$ and thus on $F_{i-1}$, so
$F_{\sigma(i+1) -1}$ depends on $F_{i-1}$. Finally, the subspace
$F_{\sigma(i+1)}$ (which in $K_A$ depended on $F_i$) does not depend
in $K_{u_iA}$ on $F_i$.

\item Note that the above two arguments are identical if exactly one
of $\sigma(i)$ and $\sigma(i+1)$ is negative.

\item Finally, if both $\sigma(i)$ and $\sigma(i+1)$ are negative,
then in the original cup diagram $CupD(A)$, both $i$ and $i+1$ are
orphans. This is exactly the case where there is no $W$-graph neighbor
to $A$ having $i$ as a descent (see Westbury\cite{We}). So we show
that the intersection $K_A \cap K_B$ has to be empty. In the component
$K_B$ where $i$ is a descent, $F_{i+1} =N^{- 1}(F_{i-1})$ and so
$F_{i+1}$ is chosen to contain exactly one lower image than $F_{i-1}$
as well as $F_i$.

On the other hand, consider a flag $F$ in the component $K_A$ where
$i$ is not a descent. In the original cup diagram $CupD(A)$, both $i$
and $i+1$ are orphans. Therefore $F_{i-1}$ is also dependent
(otherwise $i-1$ would have to connect to $i$). So $F_{i}$ is chosen
to contain exactly one lower image than $F_{i-1}$, and in fact $F_{i}
= N^{-k}(\im N^{j})$ for some $j,k$. Then $F_{i+1}$ is chosen to
contain exactly one lower image than $F_{i}$ and in fact $F_{i} = N^{-
k}(\im N^{j-1})$. Thus, in $K_A$, the subspace $F_{i+1}$ must contain
two lower images than $F_{i-1}$. But in $K_B$, the subspace $F_{i+1}$
must contain exactly one lower image then $F_{i- 1}$. Therefore there
are no flags in the intersection $K_A \cap K_B$.
\end{enumerate}

The following diagram illustrates the possibilities; in each case, the
 relevant subdiagrams of the cup diagrams $CupD(A)$ and $CupD(u_i A)$
 are on the top left and top right respectively. The relevant
 subdiagram of $CupD(B)$ is on the bottom in all cases.  \vskip12pt
 \centerline{\epsfbox{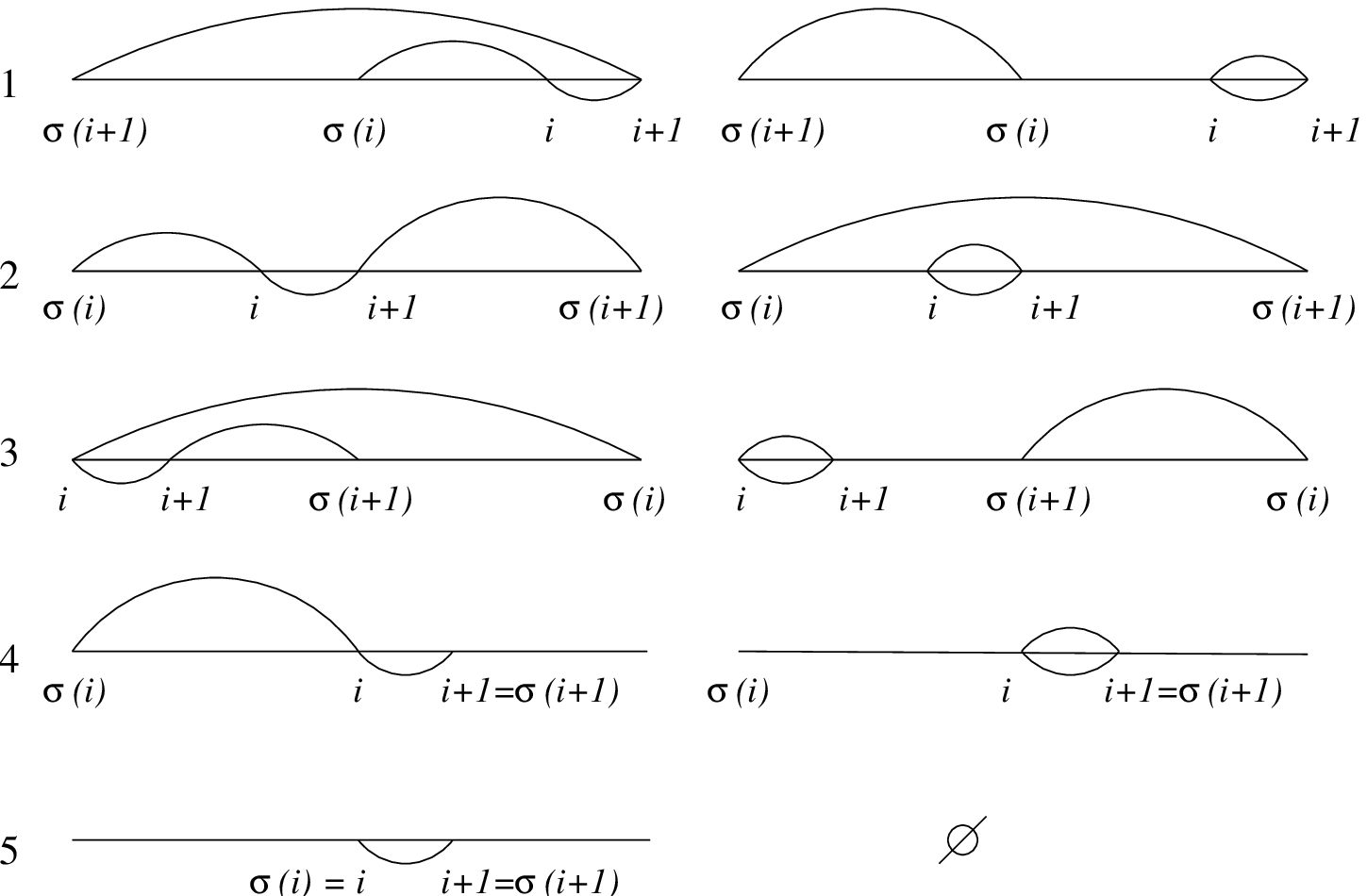}} \vskip12pt

Now we have demonstrated in all cases that the $i$-saturation of $K_A
\cap K_B$ is indeed equal to $K_{u_i A} \cap K_{B}$. Since $K_A$ is
not a union of lines of type $i$, the conclusion of Theorem
\ref{tworowchar} ensures that $K_A$ has no subvariety that consists of
lines of type $i$. Thus the Decomposition Theorem  yields
the conclusion of the theorem.
\end{proof}

As a complement, we describe the computation of the inner product
matrix of the Kazhdan-Lusztig basis for a two-row shape.

\begin{theorem} (Westbury\cite{We}, Graham-Lehrer\cite{GL}) Let $\tau$
be a two-row Young shape, and let $A$ and $B$ be two standard tableaux
on the shape of $\tau$. Consider the diagram formed by placing the cup
diagram $CupD(A)$ above the horizontal line, and $CupD(B)$
below. Suppose the diagram contains $r$ closed loops, and also that
the endpoints of each open arc are pointing in opposite
directions. Then the inner product of two Kazhdan-Lusztig basis
vectors $c_A$ and $c_B$ is $\la c_A, c_B \ra = (t + t^{-1})^r$. If an
open arc in the diagram has both ends pointing in the same direction,
then the inner product is $0$.
\end{theorem}

\begin{proof} See Westbury\cite{We}, Sections 5 and 7, and
Graham-Lehrer\cite{GL}, Section 6.  Note that their answer differs
from ours by a sign because they use the other Kazhdan-Lusztig basis
arising from the elements $C_w$.
\end{proof}

\begin{conjecture} Based on the strong evidence of the above
calculations, we conjecture that the pairwise intersections of
components of Springer fibers of two-row type are also iterated
$\mathbb{CP}^1$ bundles. It would suffice to show that each pairwise
intersection admits a description of the same form as Theorem
\ref{tworowchar}. Of course, Theorem \ref{tworowchar} shows that the
intersection of two components consists of all flags that satisfy the
descriptions of both components simultaneously. So it remains to show
that there is a single description of the same form for the
intersection.
\end{conjecture}

\section{Further speculations}

Much research into the relation between the Kazhdan-Lusztig basis and the Springer fibers for 
$GL_n(\mathbb{C})$
has been stimulated by the conjecture in Kazhdan-Lusztig\cite{KL} $6.3$, which states that with the tableau 
labelings
of the Springer fiber components and the Kazhdan-Lusztig basis vectors, the codimension $1$ 
intersections
of the components yield the edges of the  left cell $W$-graphs. It was also conceivable that the 
Springer
basis and the Kazhdan-Lusztig basis were the same at the level of $S_n$, and that perhaps
there was a way to get an Iwahori-Hecke algebra action on the Springer basis\cite{Spr2}.

Recent work of Kashiwara-Saito\cite{KS} disproved a conjecture concerning irreducibility of 
a certain characteristic variety, which implies that the Springer and 
Kazhdan-Lusztig bases are indeed different at the $S_n$ level, 
and disproves conjecture $6.3$ in general.

For hook shapes, the Kazhdan-Lusztig basis is known to coincide with the Springer basis for 
$S_n$(see G\"uemes\cite{Gu});
the conjecture \cite{KL} $6.3$ also holds for two-row shapes, because of the equations established in 
Chapter \ref{tworowschap} and the work of Lascoux-Sch\"utzenberger\cite{LS} (see also 
Westbury\cite{We} and Wolper\cite{Wo}). In fact, for representations of $S_n$ labeled by
 hooks and two-rows (and all left cells for which the
Bruhat order coincides with the weak Bruhat order), \cite{KL} $6.3$ is known in the sense that Hotta's
transformation formula(see Hotta\cite{Ho} and Borho-Brylinski-MacPherson\cite{BBM}) for the Springer 
basis
coincides with the Kazhdan-Lusztig transformation formula (see Douglass\cite{Do}).

It is not yet clear how the results of  this dissertation fit into the framework of the
above results and counterexamples. The Kazhdan-Lusztig inner products, properly normalized, are 
always polynomials that are symmetric around $0$; that is, invariant with respect to the map $t \to t^{-1}$. So we would like them to 
correspond to a method of associating a symmetric Poincar\'e polynomial to each component  (and to each
pairwise intersection of components) of the Springer fiber. In the cases in this work, all 
components and intersections were nonsingular, so the homology Poincar\'e polynomials were 
already symmetric, once shifted appropriately. 
A natural choice for a symmetric Poincar\'e polynomial to a singular variety is the intersection 
homology Poincar\'e polynomial. In the nonsingular case, intersection homology coincides with 
ordinary homology, except for the shift. Also, intersection homology satisfies the 
crucial Decomposition Theorem\cite{BBD}. However,  it appears that the natural conjecture extending Theorems \ref{hookthm} and \ref{tworowthm} using intersection homology 
alone is not correct. The first
example of  a singular component arises in $S_6$ (see Vargas\cite{V}
and Spaltenstein\cite{Spa}). This component $X$ is 
specified in our notation by the tableau
\[ \begin{matrix}
6&4 \\ 5 & 2 \\ 3 & \\ 1& 
\end{matrix}
\] and $X$ can be given the following description. The component $X$ consists of flags $F$ such 
that $F_2 \subset \ker N$,  $F_2 \cap \im N$ has dimension at least $1$, $F_4 \cap  \ker N$ has 
dimension at least 3, $F_4 \subset N^{-1}(F_2)$, and $\im N \subset F_4$. In the Spaltenstein-
Vargas subset for this component, $F_2 \cap \im N$ is a one-dimensional space $p$ and $F_4$ is 
chosen in $N^{-1}(p)$ to contain $\im N$. Since $F_4$ and $\ker N$ are both $4$-dimensional 
subspaces of the $5$-dimensional space $N^{-1}(p)$, their intersection $F_4 \cap \ker N = PL$ 
must have at least dimension $3$ and contain $\im N$. 

The singular set of  the component $X$ consists of those flags such that $F_2 = \im N$ and $F_4 = 
\ker N$. In those cases, the one-dimensional subspace  $p$ is no longer uniquely determined, nor 
is the $3$-dimensional space $PL$. If, for each flag in $X$, we choose a $1$-dimensional subspace 
in $\im N \cap F_2$ and a $3$-dimensional subspace in $F_4 \cap \ker N$, then the resulting space  
$\tilde{X}$ of such triples is a resolution of singularities of $X$. Also, the fiber over a point 
in the singular set is a $\mathbb{CP}^1 \times \mathbb{CP}^1$. The singular set has complex 
codimension $4$. Thus the resolution is semismall.

We can see that space $\tilde{X}$ has homology Poincar\'e polynomial equal to 
$[2][2][2][2][2][2][2]$ as follows. The choice of a $1$-dimensional subspace $p$ in $\im N$ is a 
$\mathbb{CP}^1$. Then the choice of a $3$-dimensional space $PL$ in $\ker N$, containing $\im N$, 
is a $\mathbb{CP}^1$. Then the choice of a $2$-dimensional space $F_2$ containing $p$ and 
contained in $PL$ is a $\mathbb{CP}^1$. The choice of a space $F_4$ containing $PL$ and contained 
in the $5$-dimensional space $N^{-1}(p)$ is another $\mathbb{CP}^1$. Finally, since $X$ is a 
union of lines of types $1$, $3$, and $5$, the other choices each contribute a $\mathbb{CP}^1$.
This exhibits $\tilde{X}$ as an iterated fiber bundle.

 The semismall Decomposition Theorem  says that the intersection homology Poincar\'e 
polynomial $IP(\tilde{X})$ of the resolution $\tilde{X}$  equals the sum of $IP$'s of strata of 
$X$, each with multiplicity equal to the number of components of the fiber over the point in the 
stratum. So we can compute $IP(X)$, since we know that the fiber over each 
point in the stratum $\mathbb{CP}^1 \times  
\mathbb{CP}^1 \times \mathbb{CP}^1$ is 
$\mathbb{CP}^1 \times \mathbb{CP}^1$. 
Thus $[2][2][2][2][2][2][2] = IP(X) + [2][2][2]$.   However, the entry in our normalized inner 
product matrix is $[2][2][2][2][2][2][2]$, which is larger the homology Poincar\'e polynomial of 
$\tilde{X}$.

However, this suggests that perhaps the inner products correspond to some other semisimple 
perverse sheaves (see Beilinson-Bernstein-Deligne\cite{BBD}) in the intersections $K_A \cap K_B$, 
since semisimple perverse sheaves also satisfy the Decomposition Theorem and have symmetric 
Poincar\'e polynomials. We also have examples of inner products (from the same shape as the above 
example) that are the sum of intersection homology Poincar\'e polynomials of multiple irreducible 
components of the corresponding intersection of two components of the Springer fiber. This lends 
further weight to the idea of using semisimple perverse sheaves on the intersection $K_A \cap 
K_B$ of two components of the Springer fiber, since the appropriate  Poincar\'e polynomial is 
obtained by summing the Poincar\'e polynomials of the irreducible components of the intersection 
$K_A \cap K_B$. The data also suggests that it would be worthwhile to investigate the structure 
of resolutions of singularities of components of Springer fibers. 

There is now $W$-graph data
available up to $S_{15}$ (see Ochiai-Kako\cite{OK}) and it would be interesting to compute the 
inner products of the 
Kazhdan-Lusztig basis vectors from them. For instance, one could check whether they satisfy the Hard Lefschetz 
theorem.

There is of course more to be done on 
the computation of the topology of the components of the Springer fibers. 
The techniques exposed here exploit the relative simplicity of the structures of the nilpotent maps
for hook and two-row types. It would be interesting to understand these components better.
Even in the two-row case, it would be worthwhile to gain more information on the structure of the
fiber bundles, for instance extensions of Lorist's theorem\cite{Lo} concerning the $e$-invariants
of the nontrivial $\mathbb{CP}^1$ bundles.

We believe that further study of the Kazhdan-Lusztig inner products and their relation to the components of 
Springer fibers  will prove to be fruitful, and that there are many questions left to be answered 
here.

\bibliographystyle{amsplain}
\bibliography{paper}

\end{document}